\documentclass[11pt]{article}

\RequirePackage{amsthm,amsmath,amsfonts,amssymb,graphicx,mathtools,enumerate,placeins}
\usepackage[round]{natbib}

\usepackage[CJKbookmarks=true,
            bookmarksnumbered=true,  
			bookmarksopen=true,
			colorlinks=true,
			citecolor=blue,
			linkcolor=blue,
			anchorcolor=red,
			urlcolor=blue]{hyperref}

\usepackage[letterpaper, left=1.2truein, right=1.2truein, top = 1.2truein, bottom = 1.2truein]{geometry}

\usepackage[blocks, affil-it]{authblk}

\usepackage[dvipsnames]{xcolor}
\usepackage{tikz}

\theoremstyle{plain}

\newtheorem{theorem}{Theorem}[section]
\newtheorem{lemma}[theorem]{Lemma}

\newtheorem{proposition}{Proposition}[section]
\newtheorem{thm}{Theorem}[section]

\newtheorem{remark}{Remark}[section]

\theoremstyle{remark}

\def\H{{\mathcal H}}

\def\Z{{\mathcal Z}}

\def\q{\alpha}

\def\iid{\textnormal{iid}}
\def\mFDR{\textnormal{mFDR}}
\def\pFDR{\textnormal{pFDR}}
\def\mlfdr{\textnormal{lfdr}}
\newcommand{\clfdr}{\textnormal{clfdr}}
\def\BH{\textnormal{BH}}

\def\PI{\textnormal{PE}}
\def\PE{\textnormal{PE}}
\def\FDP{\textnormal{FDP}}
\def\hFDP{\widehat{\textnormal{FDP}}}
\def\bFDR{\textnormal{bFDR}}
\def\FDR{\textnormal{FDR}}
\def\maxlfdr{\textnormal{max-lfdr}}
\def\rank{\textnormal{rank}}
\def\lfdr{\textnormal{lfdr}}

\def\th{\textnormal{th}}
\def\de{d}
\def\d{\delta}
\def\e{\varepsilon}
\def\q{\alpha}
\def\var{\textnormal{Var}}
\def\p{\textbf{p}}
\def\R{\mathbb{R}}
\def\E{\mathbb{E}}
\def\P{\mathbb{P}}
\def\e{\varepsilon}
\def\L{J}

\title{A frequentist local false discovery rate}
\author[1]{Daniel Xiang}
\author[1]{Jake A. Soloff}
\author[2]{William Fithian}
\affil[1]{Department of Statistics, University of Chicago} \affil[2]{Department of Statistics, University of California, Berkeley} 

\begin{document}
\maketitle

\begin{abstract}
    The local false discovery rate (lfdr) of \citet{efron2001empirical} enjoys major conceptual and decision-theoretic advantages over the false discovery rate (FDR) as an error criterion in multiple testing, but is only well-defined in Bayesian models where the truth status of each null hypothesis is random. We define a frequentist counterpart to the lfdr based on the relative frequency of nulls at each point in the sample space. The frequentist lfdr is defined without reference to any prior, but preserves several important properties of the Bayesian lfdr: For continuous test statistics, $\lfdr(t)$ gives the probability, conditional on observing {\em some} statistic equal to~$t$, that the corresponding null hypothesis is true.  Evaluating the lfdr at an individual test statistic also yields a calibrated forecast of whether its null hypothesis is true. Finally, thresholding the lfdr at $\frac{1}{1+\lambda}$ gives the best separable rejection rule under the weighted classification loss where Type I errors are $\lambda$ times as costly as Type II errors. The lfdr can be estimated efficiently using parametric or non-parametric methods, and a closely related error criterion can be provably controlled in finite samples under independence assumptions. Whereas the FDR measures the average quality of all discoveries in a given rejection region, our lfdr measures how the quality of discoveries varies across the rejection region, allowing for a more fine-grained analysis without requiring the introduction of a prior. 
\end{abstract}

\section{Introduction}\label{sec-intro}

Suppose that we are testing a scientific hypothesis, and observe a $z$-statistic equal to~$3$. How confidently can we reject the corresponding null hypothesis in favor of the alternative? This simple and natural question could hardly be better crafted to embarrass frequentist statisticians. Notwithstanding the common lay misinterpretation of the $p$-value (in this case roughly $0.0027$) as the posterior probability that the null is true in light of the data, calculating this probability in fact requires further information, namely the prior probability that the null is true and the distribution of the test statistic under the alternative. Bayesians are willing to supply these quantities, but face other difficulties: different observers' subjective beliefs may vary widely, and many scientists resist granting that the truth or falsehood of a concrete scientific hypothesis is a random event whose probability rises and falls according to an observer's prejudices \citep{goodman1999toward, savage1972foundations}.

Both frequentists and Bayesians are better equipped to answer the question when the hypothesis is one of many under consideration, provided that the other hypotheses are considered {\em relevant}, meaning informally that the cases are sufficiently alike to justify a combined analysis. Then, hierarchical or empirical Bayesian methods are appealing because they allow the subjective prior to be replaced with one that is wholly or partly learned from the data. However, calculating posterior probabilities still requires the analyst to model the truth status of individual hypotheses as random variables, and to mathematically formalize the assumption of relevance, typically by assuming that the hypotheses and test statistics are exchangeable across the cases, or else by introducing a parametric model for their dependence. In many scientific contexts, these may be difficult assumptions to accept, even if we are willing in principle to proceed under a Bayesian framework.

\begin{figure*}[t!]
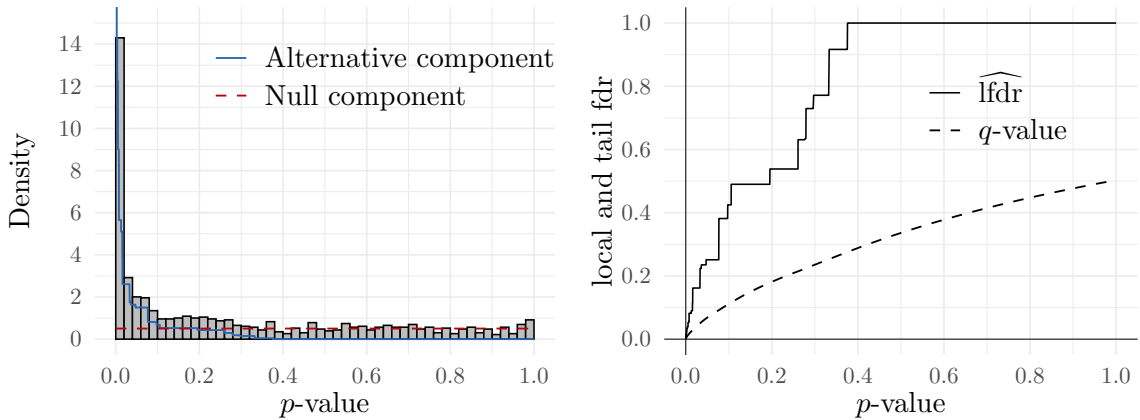

    \centering
    \begin{minipage}{.5\textwidth}
        \include{figures/1a}
    \end{minipage}%
    \begin{minipage}{.5\textwidth}
        \include{figures/1b}
    \end{minipage}
    \vspace*{-30pt}
    \caption{Microbiome preservation example. The left panel shows a histogram of permutation $p$-values comparing relative abundance in fresh vs. eight-week-old preserved samples in each of $m=1147$ species, along with a nonparametric estimate of the null and alternative components of the mixture density using the empirical Bayes estimator of \citet{strimmer2008unified}. The right panel shows the corresponding lfdr estimates, as well as the $q$-value of \citet{storey2002direct}. While these lfdr estimates would be difficult to justify as posterior probabilities in a fully Bayesian analysis, they have natural interpretations in our frequentist framework.}
    \label{fig:microbiome}
\end{figure*}

Figure~\ref{fig:microbiome} shows a histogram of 1147 permutation $p$-values from a microbiome experiment by \citet{song2016preservation}, discussed in Section~\ref{sec:microbiome-example}. Each $p$-value tests whether the relative abundance of a given bacterial species is altered by laboratory storage of biological samples, and is calculated under nonparametric assumptions justified by the controlled experimental conditions. However, the complex taxonomic and ecological relationships between the different species make it highly implausible that the $p$-values for different species are stochastically independent, or that the scientific hypotheses are exchangeable {\em a priori}, and the need to model these dependencies presents a forbidding obstacle to a full Bayesian analysis.

Frequentists reluctant to embrace Bayesian assumptions commonly choose instead to estimate or control the false discovery rate (FDR), defined as the expected fraction of true nulls among hypotheses rejected by a multiple testing method such as the Benjamini--Hochberg (BH) procedure \citep{benjamini1995controlling}. The $q$-values of \cite{storey2002direct} give a kind of FDR estimate for individual hypotheses, but they answer quite a different question than our original one: roughly, the $q$-value for a hypothesis with test statistic $t$ estimates the proportion of true nulls among all hypotheses with test statistics as extreme {\em or more extreme} than $t$. As Figure~\ref{fig:microbiome} illustrates, misinterpreting the $q$-value as a measure of confidence in a given discovery would make us systematically, and often severely, over-optimistic; see Section~\ref{sec-lfdr-and-FDR} for further discussion of this point.

Inspired by the FDR, \citet{efron2001empirical} reinterpreted the BH procedure as an empirical Bayes method, and presented nonparametric methods for estimating the (Bayesian) \emph{local false discovery rate} (lfdr), which they defined as the probability that an individual null hypothesis is true in light of the data; see Section~\ref{sec:bayes-model}. Subsequent work has developed practical methods for estimating the lfdr that avoid detailed modeling of dependence wherever possible, but has remained firmly within the Bayesian modeling paradigm. Thus, the Bayesian lfdr successfully answers our original motivating question, but it does so by reintroducing the prior that the FDR framework so deftly evades.

In this work we introduce a new definition of the lfdr that addresses our motivating question within a fully frequentist model. Suppose that we observe test statistics $z_1,\ldots,z_m$ for hypotheses $H_1,\ldots,H_m$, of which $m_0$ are true nulls. Let $f^{(i)}(t)$ denote the density of $z_i$, and define the (frequentist) lfdr as the relative frequency of null statistics at each point in the sample space:
\begin{equation}\label{eq:def-lfdr}
\lfdr(t) \;\coloneqq \sum_{i: H_i \text{ is true}} f^{(i)}(t) \;\bigg/\; \sum_{i=1}^m f^{(i)}(t).
\end{equation}

If $m$ is large and the dependence between test statistics is mild, $\lfdr(t)$ approximates the proportion of nulls among hypotheses whose statistics fall in a small neighborhood of $t$. 

In the common setting where all of the null statistics share the same density $f_0$, we obtain the simpler expression
\begin{equation}\label{eq:def-lfdr-simple}
\lfdr(t) = \bar{\pi}_0 f_0(t) \,\big/\, \bar{f}(t),
\end{equation}
where $\bar{\pi}_0 \coloneqq m_0/m$ is the true null proportion, and 
\[
\bar{f}(t) \coloneqq \frac{1}{m} \sum_{i=1}^m f^{(i)}(t)
\]
is the average density. The null density $f_0$ could, for example, represent the standard Gaussian distribution if the statistics are $z$-values, or the uniform distribution on $[0,1]$ if they are $p$-values. The lfdr estimates in Figure~\ref{fig:microbiome} estimate $\bar\pi_0$ using Storey's estimator \citep{storey2002direct}, and $\bar{f}(t)$ using Grenander's estimator of a monotone density \citep{grenander1956theory}. Let $\Z$ generically denote the common sample space where $z_1,\ldots,z_m$ are realized. For simplicity of exposition, we will assume throughout that the statistics are continuous, but most of our results extend to discrete sample spaces.

Readers may recognize the expression \eqref{eq:def-lfdr-simple} as nearly identical to the definition of lfdr in the \emph{Bayesian two-groups model} of \citet{efron2001empirical}, but there are key differences. Most importantly, our frequentist lfdr is defined only in terms of the marginal densities of the $m$ statistics, effectively replacing the Bayesian prior with the finite population of cases under study. Correspondingly, our definition of $\lfdr(t)$ is not a Bayesian posterior probability that $H_i$ is true given $z_i = t$, since that probability is either zero or one in our frequentist setting. Instead, the frequentist lfdr is the conditional probability, given {\em some} statistic equals $t$, that its null hypothesis is true:
\begin{equation}\label{eq:almost-bayes-cond}
\lfdr(t) = \P(H_J \text{ is true } \mid z_J = t, \text{ for some } J).
\end{equation}
The truth status of $H_J$ in~\eqref{eq:almost-bayes-cond} is random because the index $J$ is random.

In Section~\ref{sec:compare-bayes} we prove the relation \eqref{eq:almost-bayes-cond} and establish two other interpretations of our lfdr: First, we show that $\lfdr(z_1),\ldots,\lfdr(z_m)$ are calibrated forecasts for the unknown truth of $H_1,\ldots,H_m$. 
Second, we show that thresholding the lfdr at $1/(1+\lambda)$ implements the optimal decision rule of the form $\delta:\; \Z \to \{\text{accept, reject}\}$ for testing $H_1,\ldots,H_m$ when Type I errors are $\lambda$ times as costly as Type II errors, coinciding with the Bayes decision rule in the two-groups model. 

Section~\ref{sec-lfdr-and-FDR} compares and contrasts the lfdr with the FDR and two related criteria for frequentist error control, the marginal FDR (mFDR) and positive FDR (pFDR). While these criteria are commonly evaluated {\em globally} on the full rejection set, the characterizations of $\lfdr$ in \eqref{eq:def-lfdr} and \eqref{eq:almost-bayes-cond} can be understood as limiting forms of the mFDR and pFDR respectively, for a hypothetical {\em local} rejection rule that ``rejects'' only statistics in a small neighborhood of $t$.

Although the frequentist lfdr depends on unknown quantities, it can be estimated efficiently from the data if $m$ is reasonably large and the dependence between test statistics is not too strong, typically using the same methods developed for estimating the Bayesian lfdr. In the formulation \eqref{eq:def-lfdr-simple}, $f_0$ is typically known and $\bar{\pi}_0$ can be conservatively bounded above by $1$ or estimated using standard techniques, leaving only the problem of estimating the average density $\bar{f}$. Section~\ref{sec-estimation} discusses approaches to this problem based on standard parametric or nonparametric methods for density estimation in the i.i.d. setting, and argues that classical empirical Bayes methods commonly applied under the two-groups model can be understood as estimates of the frequentist lfdr.

Section~\ref{sec-boundary-FDR} introduces the {\em boundary FDR}, a frequentist error criterion for multiple testing defined as the probability that the last rejection (i.e., the rejection with the largest $p$-value) is a false discovery. We show that the support line procedure of \citet{soloff2024edge}, which is closely related to Strimmer's monotone lfdr estimator \citep{strimmer2008unified}, controls this criterion in finite samples under an independence assumption. Section~\ref{sec-examples} illustrates the application of the lfdr in simulation and real data examples, and Section~\ref{sec-discussion} concludes.

\section{Interpreting the frequentist lfdr}
\label{sec:compare-bayes}

\subsection{Review of the Bayes two-groups model}\label{sec:bayes-model}

In prior work, the local false discovery rate has been defined with respect to Bayesian models. The most well-known of these is the so-called {\em Bayesian two-groups model} of \citet{efron2001empirical}. In that model, each hypothesis has an independent chance $\pi_0$ of being true, and the statistic $z_i$ is distributed according to density $f_0(t)$ if $H_i$ is true, and $f_1(t)$ otherwise.

The Bayesian local false discovery rate is defined as the posterior probability that $H_i$ is true in light of the data:
\begin{align}
\label{eq:def-Bayesian-lfdr}
    \lfdr^*(t) \coloneqq \P(H_i \text{ is true} \mid z_i=t)
    = \pi_0 f_0(t)\,\big/\, f(t),   
\end{align}
where $f(t) = \pi_0 f_0(t)+(1-\pi_0) f_1(t)$ is the marginal density of $z_1,\ldots,z_m$. To avoid confusion, we use an asterisk to distinguish the Bayesian lfdr from our frequentist lfdr.

When $\pi_0, f_0,$ and $f_1$ are fixed and known, the posterior probabilities $\lfdr^*(z_1), \ldots,\lfdr^*(z_m)$ fully describe the posterior, and they represent the sharpest calibrated forecasts for the truth status of the hypotheses $H_1,\ldots,H_m$ \citep{dawid1982well,gupta2020distribution}. They also have a natural interpretation in decision theory. For a decision rule $\delta$ that returns an accept/reject decision for each hypothesis, define the {\em weighted classification loss}, which penalizes the analyst $\lambda$ for false positives and false negatives at different rates:
\begin{align}
\label{def:weighted-classification-loss}
    L_\lambda(H, \delta) \coloneqq \lambda \cdot (\# \text{false positives}) + (\# \text{false negatives}),
\end{align}
where a false positive occurs when $\delta$ rejects a true $H_i$, and a false negative occurs when $\delta$ accepts a false $H_i$. A straightforward calculation shows that the risk-minimizing rule is to reject $H_i$ if and only if $\lfdr^*(z_i)$ is below $1/(1+\lambda)$ \citep{sun2007oracle}.

This rejection rule is especially simple to interpret when the statistics are $p$-values with uniform null density $f_0(t) \equiv 1$ on $\Z = [0,1]$, in which case this rule equates to rejecting $H_i$ whenever $z_i$ is observed in a region with density $f(t) \geq (1+\lambda) \pi_0$. For example, if a false positive is $\lambda = 4$ times as costly as a false negative, then we should reject when $\lfdr^*(z_i) \leq 0.2$, or equivalently when $f(z_i) \geq 5 \cdot \pi_0$. 

The intimate connection between the Bayesian lfdr and the marginal density $f$ has a very convenient consequence in the empirical Bayes setting where $f_1$ and $\pi_0$ are unknown. If $\pi_0 \approx 1$, then estimating the marginal density $f(t)$ from the i.i.d. sample $z_1,\ldots,z_m$ is nearly equivalent to estimating $\lfdr^*(t)$, and determining the optimal rejection rule amounts to finding a super-level set of $f$.
 
In the two-groups model, these interpretations of the Bayesian lfdr are all easy consequences of standard Bayesian calculations. None of them carry over directly to our frequentist model: if the truth status of $H_1,\ldots,H_m$ is fixed, then (i) the probability $H_i$ is true in light of the data is always $0$ or $1$, (ii) the optimal forecasting rule is to forecast that the true hypotheses are true and the false ones are false, and (iii) the best decision rule for any $\lambda$ is to reject the false hypotheses and accept the true ones. Nevertheless, all three properties of the Bayesian lfdr have close analogs in our frequentist model, as we explore in the next section.

\subsection{Three interpretations of the frequentist lfdr}\label{sec:three-interp}

Section~\ref{sec-intro} gave three interpretations of the frequentist lfdr that are close analogs of properties enjoyed by the Bayesian lfdr. We now review and elaborate on them:

\vspace{+.5em}\noindent\textbf{Interpretation 1: Conditional probability.} For a fixed value $t \in \Z$, $\lfdr(t)$ is the conditional probability that a hypothesis with test statistic equal to $t$ is a true null.

\begin{thm}
\label{thm:conditional-probability}
Suppose $z_1,\dots,z_m$ are jointly absolutely continuous. Then
\begin{align*}
    \mlfdr(t) = \P(H_J \textnormal{ is true} \mid z_J=t, \textnormal{ for some } J).
\end{align*}
where $J$ is the (random) index of the statistic with $z_J=t$.
\end{thm}

When $\Z$ is discrete (in which case $f_0$ and $f$ represent probability mass functions) this property does not generalize directly in the way we might initially expect: conditional on the event that {\em at least one} index $J$ has $z_J = t$, the probability that a randomly selected one is truly null is not in general equal to $\lfdr(t)$, unless we weight the probabilities by the number of statistics equal to $t$. Instead, we have
\[
\lfdr(t) = \frac{\E[\#\{j:\; z_j = t, H_j \text{ true}\}]}{\E[\#\{j:\; z_j = t\}]}.
\]
This ratio is closely related to the marginal false discovery rate (mFDR). See Section~\ref{sec-lfdr-and-FDR} for further discussion of connections between the lfdr, FDR, and mFDR.

\vspace{+.5em} \noindent\textbf{Interpretation 2: Calibrated forecast.} The lfdr evaluated at the observed statistics $z_1,\ldots,z_m$ makes calibrated forecasts for the truth of the null hypotheses $H_1,\ldots,H_m$, where a function $g:\R \to [0,1]$ is said to be calibrated if 
\begin{align*}
\P(H_{J}\text{ is true} \mid g(z_J) = \q, \text{ for some } J) = \q.
\end{align*}
\begin{thm}
\label{thm:calibration}
Let $\ell_i \coloneqq \lfdr(z_i)$ for $i=1,\dots,m$ and suppose $z_1,\dots,z_m$ are jointly absolutely continuous. Then
\begin{align*}
    \P(H_{J} \textnormal{ is true} \mid \ell_J = \q, \textnormal{ for some } J) = \q,
\end{align*}
for any $\q \in\textnormal{range}(\mlfdr)$. Furthermore, $\mlfdr$ is the finest calibrator in the following sense: if $g: \R \to [0,1]$ is calibrated, then for any $t$,
\begin{align}
\label{eq:finest-calibrator}
    g(t) = \E ( \mlfdr(z_I) \mid g(z_I) = g(t)),
\end{align}
where $I\sim \textnormal{Uniform}\{1,\dots,m\}$.
\end{thm}

\noindent\textbf{Interpretation 3: Optimal rejection rule.} Thresholding $\lfdr(z_i)$ at $1/(1+\lambda)$ gives the optimal separable rejection rule for testing $H_1,\ldots,H_m$ under the weighted classification loss with weight $\lambda$, defined in \eqref{def:weighted-classification-loss}.

For a decision rule $\delta(z_1,\dots,z_m)\in \{\text{reject},\text{accept}\}^m$, the weighted classification risk is minimized over separable decision rules by the one that thresholds the frequentist lfdr.
\begin{thm}
\label{thm:compound-decisions}
If $\delta$ is a separable decision rule, i.e. $\delta_i(z_1,\dots,z_m) = g(z_i)$ for some univariate function $g$, then
\begin{align*}
    \E L_\lambda(H,\d) \geq \E L_\lambda(H,\mathfrak{d}^*),
\end{align*}
where 
\begin{align}
\label{eq:oracle-separable-rule}
    \mathfrak{d}_i^*(z_1,\dots,z_m) = \begin{cases}
    \textnormal{reject} \hspace{1em} &\textnormal{if } \lfdr(z_i) \leq \frac{1}{1+\lambda} \\
    \textnormal{accept} &\textnormal{otherwise}.
    \end{cases}
\end{align}    
\end{thm}

\subsection{Limitations}

The previous section discussed three interpretations of the frequentist lfdr. In practice, lfdr depends on unknown quantities such as $\bar{\pi}_0$ and $\bar{f}$ that must be estimated. Density estimation is a difficult problem in general, and therefore estimating the lfdr can be hard unless we can rely on assumptions like monotonicity or smoothness for the average density $\bar{f}$. Strong dependence between the test statistics can present an additional complication. We discuss the problem of estimating the lfdr in section~\ref{sec-estimation}, and in particular how empirical Bayes estimates target the frequentist lfdr when observations are non-i.i.d.

Another limitation of our lfdr function is that it does not account for additional information that may be known by the analyst. Formally, the conditional probability interpretation applies only to an analyst who is ignorant or indifferent about {\em which} null hypothesis corresponds to the test statistic realized at, e.g. $z_J=3$. As a result, the lfdr may not match the posterior belief of an analyst who has different prior opinions about the likelihood of different hypotheses being true. In that case, it may be more appropriate to choose a smaller reference class that represents a subset of the hypotheses under study. 

Our last interpretation assumes separability of our decision rule, but this restriction is somewhat artificial. 
In the absence of covariates, we could instead restrict to permutation equivariant (PE) rules, which implies that the rejection threshold depends only on the set of values $\{z_1,\dots,z_m\}$ and not on the order in which they are observed. In large samples with independent observations, the best PE decision rule is close to the best separable rule, mirroring a well-known phenomenon in the empirical Bayes literature (see, e.g. \cite{hannan1955asymptotic} and \cite{greenshtein2009asymptotic}). We elaborate on this point in Section~\ref{subsec-compound-lfdr} of the Appendix.

\section{lfdr and FDR}
\label{sec-lfdr-and-FDR}
\begin{figure}[tb]
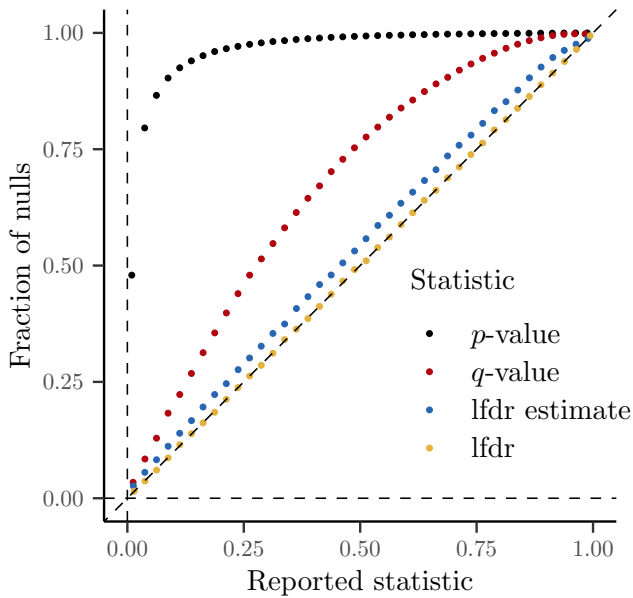

    \centering
    \include{figures/2}
    \vspace*{-30pt}
    \caption{Simulation example. The true means in \eqref{ex:gaussian-means} are $\mu_1=\dots=\mu_{m_1}=2$ and $\mu_i=0$ for $i = m_1+1,\dots,m$, where $m=3000$ and $m_1=150$, so the true null proportion is $95\%$. The calibration curve for the frequentist lfdr is plotted in yellow. We repeat the experiment $10^4$ times to assess calibration.} 
    \label{fig:ggm-calibration-plot}
\end{figure}

The lfdr was originally proposed by \citet{efron2001empirical} as a modification of the FDR, which measures the expected fraction of Type I errors among all rejections made by a multiple testing procedure. If such a procedure makes $R \geq 0$ rejections (or ``discoveries''), of which $V$ correspond to true nulls (or false discoveries), then the {\em false discovery proportion} (FDP) is defined as the realized fraction of false discoveries, and the FDR is its expectation:
\[
\FDP = \frac{V}{1 \vee R},\quad \FDR = \E [\FDP].
\]

\citet{efron2001empirical} showed that the BH procedure, the best-known method for controlling the FDR at a pre-specified level $\q$, can be interpreted as an empirical Bayes method that rejects the null for all $p$-values below a data-adaptive threshold, defined with respect to an estimator of the FDP for all rejection regions of the form $[0,t]$:
\begin{equation}\label{eq:bh-ebayes}
\hat{t}^{\text{BH}}_{\q} = \max\left\{t:\; \hFDP([0,t]) \leq \q\right\}, \quad \text{ where } \hFDP([0,t]) = \frac{mt}{\#\{i:\;p_i \leq t\}}.
\end{equation}

When the threshold $t$ is equal to an observation $p_i$, the quantity $\widehat{\FDP}([0,t])$ is a conservative estimate for the false discovery proportion among $p$-values below $p_i$.
This quantity is closely related to the $q$-value \citep{storey2004strong}, defined with respect to an FDR procedure as the smallest level $\q$ at which the procedure run at level $\q$ rejects $H_i$. For the BH procedure, the $q$-values are:
\begin{align*}
    q_{(m)} &= \widehat{\FDP}([0,p_{(m)}]) \\
    q_{(i)} &= \min\left\{\widehat{\FDP}([0,p_{(i)}]), q_{(i+1)} \right\}, \hspace{1em} i=m-1,m-2,\dots,1.
\end{align*}
Equivalently, $q_i$ is the smallest $\q \in [0,1]$ for which the BH procedure \eqref{eq:bh-ebayes} satisfies $\hat{t}_\q^{\BH} \geq p_i$. 

The $q$-value has a natural Bayesian interpretation \citep{storey2003positive} in the setting of Section \ref{sec:bayes-model}, where $q_i$ approximates the conditional probability that a null hypothesis is true, given that its $p$-value fell \textit{below} the realized value of $p_i$. As mentioned in Section \ref{sec-intro}, if we were to interpret this as the probability that $H_i$ is true in light of the data, then our confidence about individual rejections would be systematically inflated.
A more relevant quantity for assessing a null hypothesis $H_i$ is the fraction of true nulls among tests whose statistic is \textit{near} that of $H_i$, which is close to $\lfdr(p_i)$ when $m$ is large and the test statistics are independent.
To illustrate the discrepancy between these two measures of confidence, we perform a numerical experiment, generating $m=3000$ Gaussian test statistics with unit variance:
\begin{align}
\label{ex:gaussian-means}
    Z_i \sim N(\mu_i,1) \hspace{1em} \text{independently for  }i=1,\dots,m,
\end{align}
where $\mu_i \in \{0,2\}$, and suppose we want to test the null hypotheses $H_i: \mu_i=0$. From the list of $z$-statistics, we compute three summary statistics for each null hypothesis: a one-sided $p$-value, a Storey-BH $q$-value, and an lfdr estimate based on the `\texttt{fdrtool}' package \citep{strimmer2008unified}, except that we use the Storey estimator of $\bar{\pi}_0$ instead of the default estimate.
The summary statistics are binned into a grid of $[0,1]$ with bin size $2.5\%$, and in each of the forty bins we calculate the proportion of true nulls. 

Figure \ref{fig:ggm-calibration-plot} displays the results of the experiment. We see that the $q$-value systematically under-estimates the chance that the null hypothesis is true, albeit less extremely than the $p$-value does. For example, we may reject a null hypothesis $H_{i}$ at level $q=25\%$, when our actual credence in the null hypothesis should be around $50\%$. By contrast, the true lfdr is exactly calibrated, and the estimated lfdr is nearly so: among $t$-statistics for which the estimated lfdr is close to $25\%$, close to a quarter correspond to true null hypotheses. 

Two quantities that are closely related to the FDR are the positive FDR and the marginal FDR, defined in this setting as
\begin{align*}
\pFDR([0,t]) &= \E[\FDP([0,t]) \mid R([0,t]))>0], \\
    \mFDR([0,t]) &= \frac{\E  V([0,t])}{\E  R([0,t])},
\end{align*}
where $V([0,t]) = \# \{i: p_i \leq t, H_i \text{ is true}\}$ and $R([0,t]) = \# \{i : p_i \leq t\}$. If we shift our focus from a rejection region $[0,t]$ to a neighborhood $[s,t]$ of the sample space, the corresponding definitions of $\mFDR$ and $\pFDR$ are
\begin{align*}
    \mFDR([s,t]) = \frac{\E V([s,t])}{\E R([s,t])}, \hspace{1em} \pFDR([s,t]) = \E\left[ \frac{V([s,t])}{R([s,t])} \mathrel{\Big |} R([s,t]) > 0\right],
\end{align*}
where $R([s,t])$ and $V([s,t])$ are the number of rejections and Type I errors for a hypothetical procedure whose rejection region is $[s,t]$. Our frequentist lfdr can be understood as a limiting version of the $\pFDR$ and $\mFDR$ as $s$ approaches $t$. 
\begin{thm}
\label{pfdr-thm}
Suppose that $p_1,\dots,p_m$ have a continuous joint probability density. Then, for any $t$
\begin{align*}
    \lim_{\e \to 0} \mFDR([t-\e,t+\e]) &= \lim_{\e \to 0} \pFDR([t-\e,t+\e]) =\lfdr(t).
\end{align*}
\end{thm}

Under further regularity conditions, the lfdr is a calibrated forecast of the false discovery proportion among rejections with $p_i \approx t$.
In other words, among all forecasts $\{i:\lfdr(p_i)\approx \q\}$, the proportion of them for which the null is true converges to $\q$ as the number of tests tends to infinity. This statement is closely related to Theorem \ref{thm:calibration} which expresses the same calibration phenomenon in a different limiting sense; the proofs for both results can be found in Section \ref{sec-proofs}.

\begin{thm}
\label{thm-calibration-lfdp}
    Suppose $z_i \sim f^{(i)}$ independently for $i=1,2,\dots$, where $f^{(1)},f^{(2)},\dots$ is a sequence of continuous pdfs on the sample space $\mathcal{Z}$. Let $H_1,H_2,\dots \in \{0,1\}$ indicate the truth statuses of a sequence of null hypotheses, where $H_i=0$ implies $f^{(i)}=f_0$, the null density. For each $m =1,2,\dots$, define the local false discovery rate at $t$ among the first $m$ null hypotheses as 
    \begin{align*}
        \lfdr_m(t) &\coloneqq \bar{\pi}_{0,m}/\bar{f}_m(t), \hspace{1em} \text{ where }\bar{\pi}_{0,m} \coloneqq \# \{i \leq m: H_i=0\}/m \text{ and } \; \bar{f}_m \coloneqq \frac{1}{m}\sum_{i=1}^m f^{(i)}.
    \end{align*}
    For any positive sequence $(\e_m)$ satisfying $\e_m \to 0$ and $m_0 \e_m \to \infty$ as $m \to \infty$,
    \begin{align}
    \label{thm-calib-part1}
        \FDP\left(\{H_i : z_i \in [t,t+\e_m], i\leq m\} \right) - \lfdr_m(t) \stackrel{p}{\longrightarrow} 0,
    \end{align}
    for any $t$ with $f_0(t)>0$, where $m_0$ is the number of true nulls among the first $m$ hypotheses. Furthermore, letting $Z \sim f_0$, if $\lfdr_m(Z)$ has a pdf which is bounded away from zero at $\alpha$ as $m \to \infty$, then for any positive sequence $(\e_m)$ satisfying $\e_m \to 0$ and $m_0\e_m \to \infty$ as $m \to \infty$,
    \begin{align}
    \label{thm-calib-part2}
        \FDP\left(\{H_i : \lfdr_m(z_i) \in [\alpha,\alpha+\e_m], i\leq m\} \right) \stackrel{p}{\longrightarrow} \alpha.
    \end{align}
\end{thm}

\begin{remark}
    The conditions stated in Theorem \ref{thm-calibration-lfdp} are sufficient but not necessary for \eqref{thm-calib-part1} and \eqref{thm-calib-part2} to hold. In fact, the result could easily be extended to the case with dependent test statistics, as long as the dependence is not too strong. However for simplicity of exposition in Theorem \ref{thm-calibration-lfdp}, we only considered the result for the independent case.
\end{remark}

\section{Estimating the lfdr}
\label{sec-estimation}

Section \ref{sec-intro} expressed the frequentist lfdr as the intensity ratio $\bar{\pi}_0 f_0(t)/\bar{f}(t)$. Assuming $f_0$ is known, we may conservatively bound $\bar{\pi}_0 \leq 1$ or estimate it via e.g. Storey's method, reducing the problem of estimating lfdr to one of estimating the average density $\bar{f}$.

Closely related is the classical problem of estimating a density $f$, given i.i.d. observations $z_1,\dots,z_m \sim f$. Many parametric and nonparametric methods have been proposed to estimate $f$. Consider the maximum likelihood estimator
\begin{align}
\label{eq:empirical-M-estimate}
    \hat{f}_m \coloneqq \underset{f \in \mathcal{F}}{\operatorname{argmax}} \hspace{.4em} \frac{1}{m} \sum_{i=1}^m \log f(z_i),
\end{align}
where $\mathcal{F}$ is a set of candidate density functions. 

If $f$ is a monotone (non-increasing) function on $[0,1]$, which is a common assumption in multiple testing given a sequence of $p$-values (\cite{genovese2004stochastic}, \cite{strimmer2008unified}), the method of \cite{grenander1956theory} can be used to estimate $f$ using \eqref{eq:empirical-M-estimate} with
$\mathcal{F}$ equal to the set of non-increasing probability densities on $[0,1]$. For Gaussian test-statistics, \cite{kiefer1956consistency} chose $\mathcal{F}$ to be the set of Gaussian mixture densities:
\begin{align*}
    \mathcal{F} = \left\{ \int \phi(z-\mu) G(d\mu): G\text{ is a probability measure} \right\}.
\end{align*}
Both of these estimators are nonparametric in the sense that the set $\mathcal{F}$ of candidate density functions is infinite-dimensional.

A parametric approach was proposed by \cite{lindsey1974construction}, where $\mathcal{F}$ is a finite-dimensional exponential family, 
\begin{align*}
    f(z) = \exp\left\{ \sum_{j=0}^{J} \beta_j z^{j} \right\}.
\end{align*}
The resulting maximum likelihood estimate for $f$ is quite smooth for moderately sized $J$, e.g. $J=7$ is the default setting in the `\texttt{locfdr}' package of \cite{efron2011locfdr} which implements Lindsey's method as a sub-routine when estimating the lfdr.

When the observations are not i.i.d., there is no single element of $\mathcal{F}$ for which the objective in~\eqref{eq:empirical-M-estimate} matches the log-likelihood of the data. 
Nevertheless, the $M$-estimator $\hat{f}_m$ can still be computed from the sequence $z_1,\dots,z_m$. In the case of Gaussian observations, i.e. $f^{(i)} = N(\theta_i,1)$, \cite{zhang2009generalized} argues that it is sensible to estimate the average marginal density~$\bar{f}$ using~\eqref{eq:empirical-M-estimate}, taking $\mathcal{F}$ to be the set of Gaussian location mixture densities. 
We now restate his intuitive argument in the current setting.

For any candidate function $f \in \mathcal{F}$, the expectation of the objective in \eqref{eq:empirical-M-estimate} is:
\begin{align*}
    \E \left[\frac{1}{m} \sum_{i=1}^m \log f(z_i)\right] &= \int \frac{1}{m}\sum_{i=1}^m f^{(i)}(z) \log f(z) \textnormal{d}z \\
    &= \E_{\bar{f}} \log f(Z), 
\end{align*}
where $Z$ is a draw from the average density $\bar{f}$. Let $f^*_m$ denote the maximizer for the deterministic analog of~\eqref{eq:empirical-M-estimate}
\begin{align}
\label{def:population-M-estimate}
    f^*_m &\coloneqq \underset{f \in \mathcal{F}}{\operatorname{argmax}} \hspace{.4em} \E_{\bar{f}}\log f(Z) \\
    \nonumber
    &= \underset{f \in \mathcal{F}}{\operatorname{argmin}} \hspace{.4em} D(\bar{f} \mathrel{ \|} f),
\end{align}
where $D(g\| h)$ is the KL distance between two probability distributions with densities $g$ and~$h$. 

Under sufficient regularity conditions, the maximizer $\hat{f}_m$ will concentrate around $f^*_m$.
In fact, $\hat{f}_m$ will still be a consistent estimate of $f^*_m$ even if the observations are mildly dependent. As long as the objective in \eqref{eq:empirical-M-estimate} converges uniformly to the population-level objective in \eqref{def:population-M-estimate}, we will have 
$D(f^*_m\|\hat{f}_m) \stackrel{p}{\to} 0$ \citep[][Theorem 5.7]{van2000asymptotic}. We record this observation in the following proposition.

\begin{proposition}
\label{prop:uniform-convergence}
    Suppose
    $M_m^*(f) \coloneqq \E_{\bar{f}}\log f(Z)$ and $\widehat{M}_m(f) \coloneqq \frac{1}{m}\sum_{i=1}^m \log f(z_i)$ satisfy
    \begin{align*}
        \sup_{f \in \mathcal{F}} |M_m^*(f)-\widehat{M}_m(f)| \stackrel{p}{\to} 0,
    \end{align*}
    as $m \to \infty$ and suppose that $\bar{f}\in \mathcal{F}$. Then $D(\bar{f} \| \hat{f}_m) \stackrel{p}{\to} 0$.
\end{proposition}

In situations where $\bar{f} \not\in \mathcal{F}$, the $M$-estimator $\hat{f}_m$ doesn't target the average density $\bar{f}$, but instead targets the element of $\mathcal{F}$ that minimizes the KL distance to $\bar{f}$.
To ensure that $f^*_m=\bar{f}$, it is sufficient that each density $f^{(i)}$ belongs to some base class of densities $\mathcal{F}_0$, and then we take $\mathcal{F} = \text{conv}(\mathcal{F}_0)$. 
For example if we knew that each observation was normally distributed with variance 1, then the mixture density $\bar{f}$ is guaranteed to be in the set of Gaussian location mixtures.

Given an estimate of the lfdr, one may use it to perform multiple testing by rejecting all null hypotheses for which this estimate is small. The resulting rejection set aims to control the lfdr among rejected null hypotheses, but in general comes with no finite-sample guarantees. In the next section, we show that a particular estimate of the lfdr, based on the maximum likelihood estimator of \cite{grenander1956theory}, leads to a multiple testing procedure of this form that satisfies an exact bound on the false discovery probability of its last rejection. We propose this latter quantity as a new error criterion, called the boundary false discovery rate ($\bFDR$), that is distinct from the usual FDR and more aligned with the concept of lfdr. 

\section{Controlling the lfdr}
\label{sec-boundary-FDR}

To evaluate multiple testing procedures, it is natural to ask whether all the rejections are individually defensible, not just whether the list of all rejections is defensible as a whole. In a Bayesian model, this question can naturally be formulated in terms of the maximum {\em a posteriori} null probability over all the rejections. \citet{soloff2024edge} define the max-lfdr for a multiple testing procedure as the expectation of this maximum, thereby evaluating a procedure $\mathcal{R} = \{i: \text{reject }H_i\}$ according to its {\em least promising} rejection,
\begin{align*}
    \maxlfdr(\mathcal{R}) = \E \left[ \max_{i \in \mathcal{R}}\; \P(H_i \text{ is true} \mid p_i) \right].
\end{align*}
In a frequentist analysis under the fixed effects model, however, it is less obvious how to formalize what we mean by the ``least promising rejection.'' In particular, because the null probability for each hypothesis is either one or zero, the maximum is always one whenever we make any false rejections at all. 

Instead, we consider the truth status of the null hypothesis associated with the largest $p$-value within the rejection region. For a procedure $\mathcal{R}$ whose rejection region $[0,\hat{\tau}]$ contains the $R$ smallest $p$-values, the {\em boundary false discovery rate} ($\bFDR$) is defined as the probability that at least one rejection was made and the null hypothesis associated with the largest $p$-value $\leq \hat{\tau}$ is true,
\begin{align}
\label{def-bFDR}
    \bFDR(\mathcal{R}) \coloneqq \P(H_{(R)} \text{ is true}),
\end{align}
where the notation $H_{(k)}$ means the null hypothesis corresponding to the $k$th smallest $p$-value, and $H_{(0)}\coloneqq \text{false}$ by convention. As \cite{soloff2024edge} observed, the bFDR is equal to the max-lfdr in the Bayesian two-groups model (Section \ref{sec:bayes-model}) with decreasing $f_1$. \cite{genovese2002operating} analyzed a multiple testing risk function (different from our weighted classification risk), conditionally on the order statistics $(p_{(1)},\dots,p_{(m)})$. This perspective gives rise to a random permutation of $(H_1,\dots,H_m)$, denoted by $(H_{(1)},\dots,H_{(m)})$, which shifts the setting to a Bayesian one, since each $H_{(i)}$ is a non-deterministic random variable. Our boundary FDR criterion focuses on the particular hypothesis $H_{(R)}$ in this re-ordering, thereby measuring the probability that our least promising rejection is a false discovery. 

\subsection{Comparison with FDR}

\begin{figure}[t]
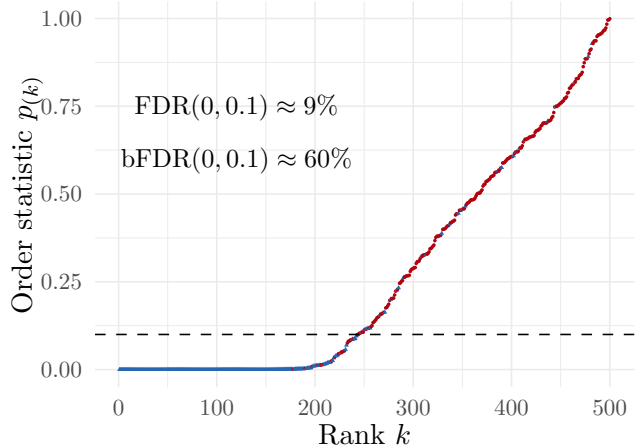

    \centering
    \include{figures/3}
    \vspace*{-30pt}
    \caption{The order statistics of $m=500$ $p$-values are plotted against their rank. They are generated with $\bar{\pi}_0=0.5$ where null $p$-values (red) are i.i.d. Uniform$(0,1)$ and alternative $p$-values (blue) are i.i.d. Beta$(0.05,1)$. The bFDR is approximated by the fraction of nulls among the largest 15 $p$-values below $0.1$.}
    \label{fig:hockey-plot}
\end{figure}

The usual FDR measures the null probability of a uniformly selected rejection:
\begin{align*}
    \FDR(\mathcal{R}) &= \P(H_{(I)} \text{ is true}),\hspace{1em} I \sim \textnormal{Uniform}\{1,\dots,R\}.
\end{align*}
Figure \ref{fig:hockey-plot} illustrates a numerical example in which the non-null $p$-values are highly concentrated near zero, leading to a substantial difference between the average-case rejection (FDR), and the ones near the boundary (bFDR).

Under a monotonicity assumption, the boundary rejection has the greatest null probability of any rejection, which implies the boundary FDR is larger than the FDR. While one might therefore be tempted to conclude that bFDR control is an inherently more conservative goal than FDR control, in practice this may or may not be the case, because one would use a larger threshold when controlling the bFDR than when controlling the FDR. For example, an analyst who equates $\lambda=4$ type II errors with a single type I error would want to control bFDR at level $1/(1+\lambda)=0.2$. The same analyst would {\em not} be satisfied with a method whose FDR is $0.2$, since the cost of the false discoveries would on average exactly cancel out the benefits of the true discoveries.

To illustrate this point, consider the weighted classification risk, which can be redefined (up to additive and multiplicative constants) as
\begin{align*}
    L_\lambda(H,\d) \coloneqq \lambda V-(R-V),
\end{align*}
where $V$ is the number of false positives among the $R$ discoveries. Taking $\lambda=4$, a procedure targeting a false discovery rate 
$V/R = 1/(1+\lambda) = 0.2$ 
achieves the same loss as a trivial procedure that simply sets $V=R=0$.

Instead, such an analyst would always aim to control FDR at some level smaller than $0.2$, for example $0.1$ so that they achieve some net benefit from the experiment. As a result, no sensible analyst would ever be interested in bFDR control and FDR control at the same level. Since bFDR control and FDR control typically wouldn't be carried out at the same level, it is unclear which is more conservative in any given case.

\subsection{Controlling the boundary FDR}

\cite{soloff2024edge} proposed the Support Line (SL) method for controlling the max-lfdr under a monotonicity constraint. The procedure run at level $\q$ rejects the $R_\q$ smallest $p$-values, where
\begin{align}
\label{num-SL-rejections}
    R_\q &\coloneqq \underset{k =0,\dots,m}{\operatorname{argmax}} \hspace{.4em} \left\{ \frac{\q k}{m} - p_{(k)} \right\}, \hspace{1em} p_{(0)} \coloneqq 0.
\end{align}
The SL method controls its boundary FDR when the nulls are independent.

\begin{thm}
\label{main-RFDR-control}
If $p_1,\dots,p_m$ are independent and $p_i \sim \textnormal{Uniform}(0,1)$ when $H_i$ is true, then
\begin{align*}
    \bFDR(\mathcal{R}_{\q}) = \bar{\pi}_0 \q,
\end{align*}
where $\mathcal{R}_{\q} \coloneqq \{i: p_i \leq p_{(R_\q)}\}$ is defined by (\ref{num-SL-rejections}).
\end{thm}

\begin{proof}[Proof of Theorem \ref{main-RFDR-control}.]

The event $\{H_{(R_\q)} \text{ is true}\}$ can be written as a disjoint union,
\begin{align*}
    \{H_{(R_\q)} \text{ is true}\} &= \bigcup_{i:H_i\text{ is true}} \{p_{(R_\q)}=p_i\}
\end{align*}
which implies 
\begin{align*}
\P(H_{(R_\q)}\text{ is true}) &= \sum_{i:H_i\text{ is true}} \P(p_{(R_\q)}=p_i) = m_0 \cdot \frac{\q}{m}.
\end{align*}
The last equality follows from Lemma 2 of \cite{soloff2024edge}, which states that for any configuration of the other $p$-values $p_1,\dots,p_{m-1}$, the probability that a null $p$-value $p_m$ achieves the optimum in (\ref{num-SL-rejections}) is equal to $\frac{\q}{m}$\footnote{An alternative proof of the fact ``$p_i \sim \text{Uniform}(0,1)\Rightarrow \P(p_{(R_\q)}=p_i) = \frac{\q}{m}$'' can be found in Section \ref{sec-proofs}.}.
\end{proof}

The SL procedure can be understood as a thresholding procedure, which thresholds a plug-in estimate of lfdr. When the null distribution is Uniform$(0,1)$, $\lfdr(t)$ is upper bounded by $1/\bar{f}(t)$ since $\bar{\pi}_0\leq 1$. The largest $p$-value in the rejection region $p_{(R_\q)}$ is equal to
\begin{align*}
    \hat{\tau}_\q \coloneqq \max\left\{ t\in[0,1]:1/\hat{f}_m(t) \leq \q \right\},
\end{align*}
where $\hat{f}_m$ is the non-parametric likelihood estimator of $\bar{f}$ \citep{grenander1956theory}, defined by \eqref{eq:empirical-M-estimate} where $\mathcal{F}$ is the set of non-increasing densities on $[0,1]$; see also Section 3.2 of \cite{soloff2024edge}. 
Under regularity conditions, for large $m$,
\begin{align*}
    \bFDR(\mathcal{R}_\q)=\E\big[ \P(H_{(R_\q)}\text{ is true} \mid p_{(R_\q)} )\big] \approx \lfdr(p_{(R_\q)}).
\end{align*}
Under mild regularity conditions, the approximation above gets better as the number of tests increases, as summarized by the following result (proved in Section \ref{sec-proofs}).

\begin{thm}
\label{flfdr-control}
Suppose $p_1,\dots,p_m$ are independent, where each $f^{(i)}$ is a continuous probability density function, equal to $1_{[0,1]}$ when $H_i$ is true, and that $\bar{f}$ has a unique solution $\tau^*_\q$ to the equation $\bar{f}(\tau^*_\q)=\q^{-1}$. If $\bar{f}$ is decreasing, and for some constants $\d,\L>0$ we have $\L \leq |\bar{f}'(t)|\leq \L^{-1}$ for all $t$ with $|t-\tau^*_\q| \leq \e$, where $\e \coloneqq \left(\frac{48}{\q \L^2}\right)^{1/3} m^{-1/3}\log(2m/\d)$, then for a constant $C>0$ depending on $\q,\L$ and $\d$,
\begin{align*}
\P\left(\left|\mlfdr(\hat{\tau}_\q)-\bar{\pi}_0 \q\right| > C m^{-1/3}\log(m/\d)\right) \leq \d.
\end{align*}
\end{thm}

\subsection{Non-Uniform$(0,1)$ null distribution}

In the argument for Theorem \ref{main-RFDR-control}, we showed that the boundary FDR control of the SL procedure is controlled when each null density is bounded,
\begin{align}
    \label{bounded-null-condition}
        H_i \text{ is true} \Rightarrow f^{(i)}(t) \leq 1 \hspace{1em} \text{ for all } t \in [0,\q].
\end{align}
This condition is distinct from requiring the nulls be super-uniformly distributed, which is an assumption commonly made in the multiple testing literature and is not sufficient in general to guarantee boundary FDR control\footnote{A counterexample to the conjecture that bFDR is controlled by super-uniform $p$-values is given in Section \ref{subsubsec-CE-super-uniform}}. 
Further regularity conditions on the super-uniform null distribution restore the guarantee $\bFDR(\mathcal{R}_\q) \leq \bar{\pi}_0 \q$. For example, when the $p$-values are generated from a one-sided Gaussian location testing problem,
\begin{align*}
        X_i \sim N(\theta_i,1), \hspace{1em} i=1,\dots,m
\end{align*}
where $H_i \text{ is true} \Rightarrow \theta_i \leq 0$. In this case, the probability density function for $p_i=1-\Phi(X_i)$ satisfies \eqref{bounded-null-condition} under the null, for $\q = 1/2$. This observation extends to one-parameter exponential families with continuous densities. The proof of this proposition is recorded in Section \ref{sec-proofs}.
\begin{proposition}
\label{prop-exp-fam}
    Let $(g_{\theta})_{\theta \in \R}$ denote an exponential family of continuous distributions on $\R$ with densities
\begin{align*}
    g_\theta(z) = \exp(\theta z - A(\theta)) g_{\theta_0}(z), \hspace{1em} \theta,z \in \R,
\end{align*}
with corresponding cdfs $(G_\theta)$.
For one-sided testing of the hypotheses $H_{i}: \theta_i \leq \theta_0$, let $\q^* = 1-G_{\theta_0}(\E_{\theta_0}Z)$ be the upper percentile of the mean under $\theta_0$. Then the null density of the one-sided $p$-value $p = 1-G_{\theta_0}(Z)$ is bounded by 1 on $[0,\q^*]$, for all $\theta \leq \theta_0$.
\end{proposition}

Another common setting in which super-uniformly distributed $p$-values arise is the discrete case, where
\begin{align*}
    p_i \sim \text{Uniform}\{1/L,\dots,L/L\}\text{  when $H_i$ is true},
\end{align*}
for some large fixed grid length $L$, e.g. the number of permutations used to compute a $p$-value for a permutation test. In this case, the boundary FDR of the SL procedure is not controlled in finite samples (see Section \ref{subsubsec-CE-discrete-uniform} for a counterexample). However, holding $m$ fixed as $L \to \infty$, or holding $L$ fixed as $m\to\infty$, the bFDR of the SL procedure is asymptotically controlled below $\alpha$. These results are stated precisely in Section \ref{subsec-discrete-uniform}, and proved in Section \ref{sec-proofs}. If one is particularly concerned about retaining a finite-sample guarantee in this setting, they may run the SL procedure on perturbed $p$-values, i.e. $\tilde{p}_i \mid p_i = \ell/L \sim \text{Uniform}((\ell-1)/L,\ell/L)$, whose null distribution is exactly Uniform$(0,1)$.

\section{Applications}
\label{sec-examples}

\subsection{Example 1: Gaussian graphical model}
\label{sec-ggm-example}
The Gaussian graphical model is an example of a setting where the frequentist lfdr is useful because the Bayesian approach requires complicated modeling, and the $q$-value approach is inherently biased.

In this model, the data arrive as $n$ i.i.d. copies of
\begin{align*}
    X \sim N(0,\Omega^{-1}),
\end{align*}
where $\Omega$ is a $d \times d$ dimensional precision matrix. This matrix encodes conditional dependence relationships between the coordinates of $X$ as follows: $\Omega_{i,j}$ is zero when the $i^{\th}$ and $j^{\th}$ coordinate of $X$ are conditionally independent, given the rest of the coordinates. To decide whether or not to reject the null hypothesis:
\begin{align*}
    H_{ij} : \; \Omega_{ij}=0,\hspace{1em} i\neq j
\end{align*}
we may compute a $t$-statistic on $n-d$ degrees from the linear model obtained by regressing $X_{j}$ against $X_{-j}$, taking $t_{ij}$ to be the standardized coefficient for $X_i$ in the fitted model.

For each pair $(i,j)$, we have $\Omega_{ij}=\Omega_{ji}$ and $t_{ij}=t_{ji}$ so the total number of hypotheses is ${d \choose 2}$. It would be inappropriate to model the $t$-statistics as independent, since $t_{ij}$ and $t_{jk}$ being large and positive is informative about the value of $t_{ik}$. In general, $t_{ij}$ is not a sufficient statistic for testing $H_{ij}$, and the posterior probability of the null $H_{ij}$ could be a complicated function of the entire sample covariance matrix.

\begin{figure}[tb]
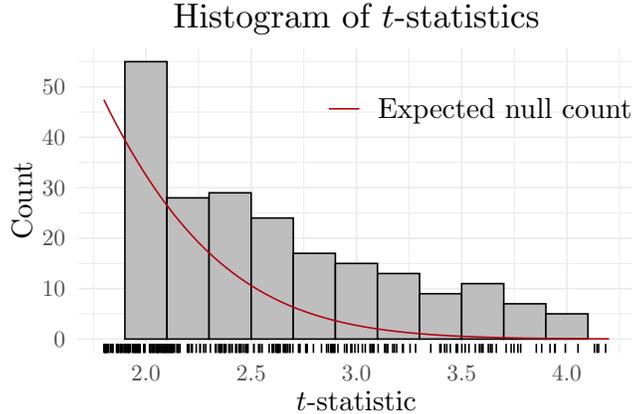

    \centering
    \include{figures/4}
    \vspace*{-30pt}
    \caption{Histogram of $t$-statistics in the GGM example, with dimensions: $d=80$, $n=10d$, $m=3,160$. The null distribution is $t_{n-d}$.}
    \label{fig:ggm-histogram}
\end{figure}

We can bypass some of these stumbling blocks by calculating a frequentist $q$-value for each $t$-statistic, but this can also be misleading. In general, the $q$-value for $t_{ij}$ substantially under-estimates the chance that $H_{ij}$ is true. 
Figure \ref{fig:ggm-histogram} shows part of the histogram of $t$-statistics generated from the previously described regression method in a Gaussian graphical model with $d=80$ and $n=10d$. 
Looking at the histogram, it is clear that we can estimate a local null proportion based on the $t$-statistics. To do so, we first calculate the expected number of null observations at, e.g. $t=3$. Overlaid in red is the Student-$t_{n-d}$ density weighted by the number of true nulls, in this case $m_0=0.95m$.
Dividing by the height of the histogram there 
yields a rough and ready estimate of the lfdr.
Compared to the BH $q$-value, which is around $6\%$ for a $t$-statistic near 3, the histogram-based estimate
of the lfdr is much higher, closer to $20\%$. 

\subsection{Example 2: Microbiome data analysis}\label{sec:microbiome-example}

This section discusses a data set from \citet{song2016preservation} on storage techniques for biological samples in microbiome analysis. In scientific investigations with microbiome data, it can be necessary to store biological samples for some period of time after collection. A key question about the integrity of the subsequent analysis is whether the relative abundance of different microbial species shifts significantly during the storage period, and whether some storage methods are better than others.

In the data set we analyze, fecal samples from six human participants were stored in 95\% ethanol solution, with microbial abundances measured using DNA sequencing techniques, both when the sample was fresh and after eight weeks of storage, with five replications of measurements per participant in each storage condition. In each replication, and for each of $9719$ microbial species\footnote{We use the same method as the original investigators for operationally defining ``species,'' which are referred to more precisely as {\em operational taxonomic units} in the microbiome literature.}, the relative abundance in each fecal sample is measured as the number of individual bacteria in that species sequenced divided by the total number of bacteria. For each species, we can use a permutation test to assess the null hypothesis that its relative abundance is unchanged after eight weeks of storage.

Due to the careful experimental design, a stratified permutation test is well-suited to test the hypothesis that the relative abundance of a given species is independent of the storage condition (fresh or eight weeks old) given the identity of the human participant. We calculate a $p$-value $p_i$ for species $i$ based on the Wilcoxon signed rank statistic, where ranks are calculated for the relative abundance of that species within each stratum. 
Because some species are sparsely observed, we restrict our analysis to the $m = 1147$ species for which the relative abundance is above zero in at least ten total replications. The permutation test is marginally valid for each species, under a generic nonparametric model. Figure~\ref{fig:microbiome} shows the $p$-value histogram as well as a nonparametric estimate for the lfdr and tail FDR, due to \citet{strimmer2008unified}. To believe in these estimates we need rely only on the assumption that the null $p$-values are approximately uniform, and that the heights of the histogram bars are informative about the mixture density (or more precisely that the empirical CDF is a good estimator for the true mixture CDF). Both assumptions appear sensible in this case without our needing to appeal to a Bayesian model.

By contrast, it would be highly challenging to specify a convincing Bayesian model for the joint distribution of $9719$ species' relative abundances under the two storage conditions. In particular, given the taxonomic structure of the different species, it is highly unlikely that the true effects of storage on each species are exchangeable, or that the observed relative abundances are independent conditional on the true effects. By shifting to our frequentist perspective, we can sidestep the difficulties of Bayesian modeling.

\subsection{Example 3: Aggregate analysis of nudges}
\label{sec-nudge}

The concept of nudging is described by \cite{thaler2009nudge} as a way of influencing people's behavior in a predictable way without restricting their options or altering economic incentives. To evaluate the overall effectiveness of psychological nudging on human behavior, \cite{mertens2022effectiveness} collected data from 447 nudge experiments in the behavioral psychology literature. 
The formulation of this question and the authors' conclusion was the subject of some debate (see e.g. \cite{maier2022no}, \cite{mertens2022reply}, \cite{szaszi2022no}).

\begin{figure}[t]
    \centering
    % Created by tikzDevice version 0.12.6 on 2024-09-28 11:49:25
% !TEX encoding = UTF-8 Unicode
\begin{tikzpicture}[x=1pt,y=1pt]
\definecolor{fillColor}{RGB}{255,255,255}
\path[use as bounding box,fill=fillColor,fill opacity=0.00] (0,0) rectangle (245.72,166.22);
\begin{scope}
\path[clip] (  0.00,  0.00) rectangle (245.72,166.22);
\definecolor{drawColor}{RGB}{255,255,255}
\definecolor{fillColor}{RGB}{255,255,255}

\path[draw=drawColor,line width= 0.9pt,line join=round,line cap=round,fill=fillColor] (  0.00,  0.00) rectangle (245.72,166.22);
\end{scope}
\begin{scope}
\path[clip] (  0.00,  0.00) rectangle (245.72,166.22);
\definecolor{fillColor}{RGB}{255,255,255}

\path[fill=fillColor] ( 25.83, 36.30) rectangle (237.22,157.72);
\definecolor{fillColor}{RGB}{190,190,190}

\path[fill=fillColor] ( 25.83, 36.30) rectangle ( 32.95,157.72);

\path[fill=fillColor] ( 32.95, 36.30) rectangle ( 40.08, 50.79);

\path[fill=fillColor] ( 40.08, 36.30) rectangle ( 47.20, 51.70);

\path[fill=fillColor] ( 47.20, 36.30) rectangle ( 54.32, 43.55);

\path[fill=fillColor] ( 54.32, 36.30) rectangle ( 61.45, 44.45);

\path[fill=fillColor] ( 61.45, 36.30) rectangle ( 68.57, 39.01);

\path[fill=fillColor] ( 68.57, 41.73) rectangle ( 75.69, 41.73);

\path[fill=fillColor] ( 75.69, 44.45) rectangle ( 82.82, 44.45);

\path[fill=fillColor] ( 82.82, 36.30) rectangle ( 92.63, 37.20);

\path[fill=fillColor] ( 92.63, 36.30) rectangle (100.24, 39.01);

\path[fill=fillColor] (100.24, 36.30) rectangle (107.85, 40.83);

\path[fill=fillColor] (107.85, 36.30) rectangle (115.46, 40.83);

\path[fill=fillColor] (115.46, 36.30) rectangle (123.07, 39.92);

\path[fill=fillColor] (123.07, 36.30) rectangle (130.68, 40.83);

\path[fill=fillColor] (130.68, 36.30) rectangle (138.29, 39.92);

\path[fill=fillColor] (138.29, 36.30) rectangle (145.90, 38.11);

\path[fill=fillColor] (145.90, 36.30) rectangle (153.51, 38.11);

\path[fill=fillColor] (153.51, 36.30) rectangle (161.12, 40.83);

\path[fill=fillColor] (161.12, 36.30) rectangle (168.73, 39.01);

\path[fill=fillColor] (168.73, 36.30) rectangle (176.34, 38.11);

\path[fill=fillColor] (176.34, 36.30) rectangle (183.95, 39.92);

\path[fill=fillColor] (183.95, 36.30) rectangle (191.56, 39.01);

\path[fill=fillColor] (191.56, 36.30) rectangle (199.17, 39.92);

\path[fill=fillColor] (199.17, 36.30) rectangle (206.78, 37.20);

\path[fill=fillColor] (206.78, 36.30) rectangle (214.39, 37.20);

\path[fill=fillColor] (214.39, 36.30) rectangle (222.00, 40.83);

\path[fill=fillColor] (222.00, 36.30) rectangle (229.61, 36.30);

\path[fill=fillColor] (229.61, 36.30) rectangle (237.22, 36.30);
\definecolor{fillColor}{RGB}{177,4,14}

\path[fill=fillColor] ( 25.83, 36.30) rectangle ( 32.95, 36.30);

\path[fill=fillColor] ( 32.95, 36.30) rectangle ( 40.08, 36.30);

\path[fill=fillColor] ( 40.08, 36.30) rectangle ( 47.20, 36.30);

\path[fill=fillColor] ( 47.20, 36.30) rectangle ( 54.32, 36.30);

\path[fill=fillColor] ( 54.32, 36.30) rectangle ( 61.45, 36.30);

\path[fill=fillColor] ( 61.45, 36.30) rectangle ( 68.57, 36.30);

\path[fill=fillColor] ( 68.57, 36.30) rectangle ( 75.69, 41.73);

\path[fill=fillColor] ( 75.69, 36.30) rectangle ( 82.82, 44.45);

\path[fill=fillColor] ( 82.82, 36.30) rectangle ( 92.63, 36.30);

\path[fill=fillColor] ( 92.63, 36.30) rectangle (100.24, 36.30);

\path[fill=fillColor] (100.24, 36.30) rectangle (107.85, 36.30);

\path[fill=fillColor] (107.85, 36.30) rectangle (115.46, 36.30);

\path[fill=fillColor] (115.46, 36.30) rectangle (123.07, 36.30);

\path[fill=fillColor] (123.07, 36.30) rectangle (130.68, 36.30);

\path[fill=fillColor] (130.68, 36.30) rectangle (138.29, 36.30);

\path[fill=fillColor] (138.29, 36.30) rectangle (145.90, 36.30);

\path[fill=fillColor] (145.90, 36.30) rectangle (153.51, 36.30);

\path[fill=fillColor] (153.51, 36.30) rectangle (161.12, 36.30);

\path[fill=fillColor] (161.12, 36.30) rectangle (168.73, 36.30);

\path[fill=fillColor] (168.73, 36.30) rectangle (176.34, 36.30);

\path[fill=fillColor] (176.34, 36.30) rectangle (183.95, 36.30);

\path[fill=fillColor] (183.95, 36.30) rectangle (191.56, 36.30);

\path[fill=fillColor] (191.56, 36.30) rectangle (199.17, 36.30);

\path[fill=fillColor] (199.17, 36.30) rectangle (206.78, 36.30);

\path[fill=fillColor] (206.78, 36.30) rectangle (214.39, 36.30);

\path[fill=fillColor] (214.39, 36.30) rectangle (222.00, 36.30);

\path[fill=fillColor] (222.00, 36.30) rectangle (229.61, 36.30);

\path[fill=fillColor] (229.61, 36.30) rectangle (237.22, 36.30);
\definecolor{drawColor}{RGB}{177,4,14}

\path[draw=drawColor,line width= 0.6pt,line join=round] ( 82.82, 36.30) -- ( 82.82,157.72);

\node[text=drawColor,anchor=base,inner sep=0pt, outer sep=0pt, scale=  0.90] at ( 85.67, 26.21) {$\tau_{q}^{\textnormal{BH}}$};

\node[text=drawColor,anchor=base,inner sep=0pt, outer sep=0pt, scale=  0.90] at (137.52, 60.39) {$\widehat{\textnormal{FDP}}([0.20,0.27])=0.32$};
\end{scope}
\begin{scope}
\path[clip] (  0.00,  0.00) rectangle (245.72,166.22);
\definecolor{drawColor}{RGB}{0,0,0}

\path[draw=drawColor,line width= 0.9pt,line join=round] ( 25.83, 36.30) --
	( 25.83,157.72);
\end{scope}
\begin{scope}
\path[clip] (  0.00,  0.00) rectangle (245.72,166.22);
\definecolor{drawColor}{RGB}{0,0,0}

\path[draw=drawColor,line width= 0.9pt,line join=round] ( 25.83, 36.30) --
	(237.22, 36.30);
\end{scope}
\begin{scope}
\path[clip] (  0.00,  0.00) rectangle (245.72,166.22);
\definecolor{drawColor}{gray}{0.20}

\path[draw=drawColor,line width= 0.9pt,line join=round] ( 25.83, 32.05) --
	( 25.83, 36.30);

\path[draw=drawColor,line width= 0.9pt,line join=round] ( 68.11, 32.05) --
	( 68.11, 36.30);

\path[draw=drawColor,line width= 0.9pt,line join=round] (110.39, 32.05) --
	(110.39, 36.30);

\path[draw=drawColor,line width= 0.9pt,line join=round] (152.66, 32.05) --
	(152.66, 36.30);

\path[draw=drawColor,line width= 0.9pt,line join=round] (194.94, 32.05) --
	(194.94, 36.30);

\path[draw=drawColor,line width= 0.9pt,line join=round] (237.22, 32.05) --
	(237.22, 36.30);
\end{scope}
\begin{scope}
\path[clip] (  0.00,  0.00) rectangle (245.72,166.22);
\definecolor{drawColor}{gray}{0.30}

\node[text=drawColor,anchor=base,inner sep=0pt, outer sep=0pt, scale=  0.80] at ( 25.83, 23.14) {0.0};

\node[text=drawColor,anchor=base,inner sep=0pt, outer sep=0pt, scale=  0.80] at ( 68.11, 23.14) {0.2};

\node[text=drawColor,anchor=base,inner sep=0pt, outer sep=0pt, scale=  0.80] at (110.39, 23.14) {0.4};

\node[text=drawColor,anchor=base,inner sep=0pt, outer sep=0pt, scale=  0.80] at (152.66, 23.14) {0.6};

\node[text=drawColor,anchor=base,inner sep=0pt, outer sep=0pt, scale=  0.80] at (194.94, 23.14) {0.8};

\node[text=drawColor,anchor=base,inner sep=0pt, outer sep=0pt, scale=  0.80] at (237.22, 23.14) {1.0};
\end{scope}
\begin{scope}
\path[clip] (  0.00,  0.00) rectangle (245.72,166.22);
\definecolor{drawColor}{RGB}{0,0,0}

\node[text=drawColor,anchor=base,inner sep=0pt, outer sep=0pt, scale=  1.00] at (131.52, 10.44) {Selection-adjusted $p$-values};
\end{scope}
\begin{scope}
\path[clip] (  0.00,  0.00) rectangle (245.72,166.22);
\definecolor{drawColor}{RGB}{0,0,0}

\node[text=drawColor,rotate= 90.00,anchor=base,inner sep=0pt, outer sep=0pt, scale=  1.00] at ( 15.39, 97.01) {Counts};
\end{scope}
\end{tikzpicture} 
    \vspace*{-30pt}
    \caption{Shown above is the histogram of one-sided $p$-values falling below $0.025$, adjusted for selection by multiplying by 40 (the reciprocal of the 2.5\% one-sided significance threshold). 
    The $\BH(q)$ threshold for $q=10\%$ 
    is around $0.27$ (or $\approx 0.007$ on the scale of the unadjusted $p$-values), below which there are 202 rejections. The estimated FDP near the edge of the rejection set (red) is around 32\%.}
    \label{fig:FDR-est-hist}
\end{figure}
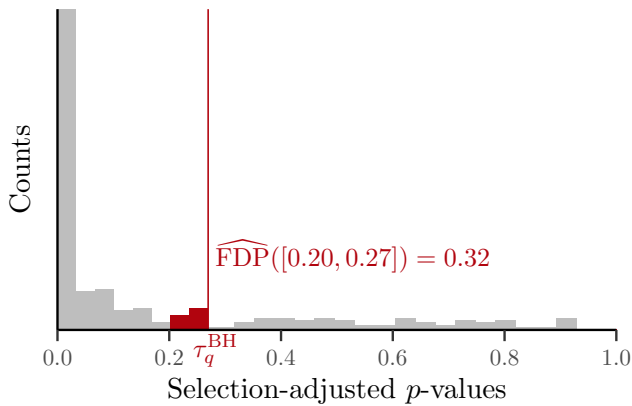

To understand the degree to which false discoveries are present in the aggregated dataset, we estimate the false discovery rate (FDR) using the Storey estimator \citep{storey2002direct} for the proportion of true nulls, restricting attention to just the $m=261$ many $p$-values falling below the $5\%$ two-sided significance level. This restriction is a way to work around the publication bias present in scientific journals; although ineffective nudges may be under-represented among published studies, the null hypotheses whose $p$-values fall within the significance region are less prone to censorship (\cite{hung2020statistical}, \cite{jaljuli2022quantifying}). 

The Storey estimator of the null proportion within the significance region is around $28\%$, suggesting that roughly a quarter of the $m=261$ results reported below the 2.5\% one-sided significance level are false discoveries. To mitigate the high rate of false claims, we ran the Storey-adjusted BH procedure \citep{storey2004strong} targeting a $10\%$ FDR, yielding a more stringent rejection threshold, as shown in Figure~\ref{fig:FDR-est-hist}, below which there are only 202 $p$-values. 

\begin{figure}[t]
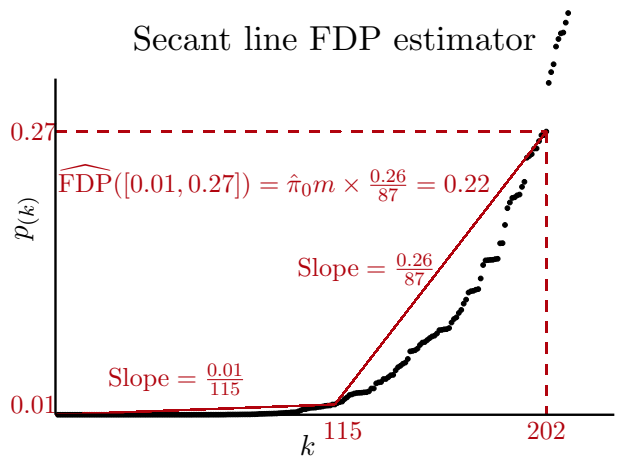

    \centering
    \include{figures/6}
    \vspace*{-30pt}
    \caption{Nudge example. The order statistics of the selection-adjusted $p$-values are plotted against their rank. 
    The estimated FDP among the smallest 202 $p$-values is close to $10\%$. Within the first half, the estimate is $0.6\%$, whereas in the second half it is $22\%$. 
    }
    \label{fig:FDR-est-secant}
\end{figure}

Upon inspecting the histogram left of the BH threshold, we find that the the estimated rate of false discoveries varies substantially. 
To estimate the rate of false discoveries within an interval $[s,t]$, we multiplied the total estimated number of nulls by the length of the interval $t-s$, and divided by the total number of observations between $s$ and $t$. This is a sensible estimate for the mFDR in $[s,t]$ because if the null distribution is Uniform$(0,1)$, then $\hat{m}_0 (t-s)$ estimates the expected number of false discoveries with $p$-values between $s$ and $t$. Visually, this estimate is proportional to the slope of a secant line drawn over the interval $[s,t]$, as illustrated in Figure~\ref{fig:FDR-est-secant} for the nudge data.
As Figure~\ref{fig:FDR-est-hist} shows, the estimated proportion of false discoveries (FDP) grossly exceeds 10\% for a subset of rejections near the rejection threshold. 

\section{Discussion}\label{sec-discussion}

\citet{efron2019bayes} remarked that ``considering the enormous gains potentially available from empirical Bayes methods, the effects on statistical practice have been somewhat underwhelming." One barrier to the wider adoption of empirical Bayes is its philosophical status. 
Frequentists have legitimate concerns about the Bayesian side of empirical Bayes.
This paper introduces a frequentist counterpart to the local false discovery rate that addresses these concerns while preserving key properties that make the Bayesian lfdr appealing: it is firmly rooted in decision theory, provides interpretable probability statements about individual hypotheses, and can be estimated using standard empirical Bayes techniques.

We close our discussion by highlighting some promising directions for future work.

\begin{itemize}
    \item \textbf{Frequentist posteriors.} The Bayesian local false discovery rate simply characterizes the posterior distribution of a binary latent variable, $H_i$, and our frequentist definition corresponds to the oracle Bayes posterior of \cite{efron2019bayes}. 
    It may be of interest to estimate the full oracle Bayes posterior beyond binary settings.
    In the Gaussian sequence model, compound decision theory has mostly focused on estimating the mean of the posterior \citep{zhang2009generalized,jiang2009general}. \vspace{.5em}
    \item \textbf{Estimation in the frequentist model.} 
    While we give one asymptotic result on estimating the lfdr (Proposition~\ref{prop:uniform-convergence}), finite-sample estimation error is a serious concern. 
    When the statistics are independent but not identically distributed, the empirical distribution is in a strong sense less dispersed than its i.i.d. counterpart \citep[see, e.g.,][Chapter~25]{shorack2009empirical}.
    It would be interesting to investigate whether this observation allows us to translate empirical Bayes guarantees into compound decision theory guarantees where the i.i.d. assumption is violated \citep{hannan1955asymptotic,han2024approximate}. 
\end{itemize}

\FloatBarrier

\section*{Acknowledgements}

This project was originally inspired by a question posed by Brad Efron after a seminar in 2018: what, if anything, does the lfdr mean if we do not believe in our Bayesian assumptions? 
A frequentist lfdr also appeared in an earlier version of \cite{weinstein2019nonparametric}, which inspired a continued development of the frequentist interpretation of lfdr.
We thank Bradley Efron, Rina Foygel Barber, Chao Gao, Ruth Heller, Peter McCullagh, Etienne Roquain, Asaf Weinstein, Stefan Wager, and Dani Yekutieli for helpful comments, discussions and support. 
J. A. Soloff gratefully acknowledges the support of the National Science Foundation via grant DMS-2023109, the Office of Naval Research via grant N00014-20-1-2337, and the Margot and Tom Pritzker Foundation.

\vspace{.5em}

\noindent\textbf{Reproducibility.} Code to reproduce all figures is available at our \href{https://github.com/dan-xiang/dan-xiang.github.io/tree/master/frequentist-lfdr-paper}{Github repository}\footnote{full link: https://github.com/dan-xiang/dan-xiang.github.io/tree/master/frequentist-lfdr-paper}.

\appendix
\section{Appendix}
\label{appn}

\subsection{Compound lfdr}
\label{subsec-compound-lfdr}

In this section, we refer to formula \eqref{eq:def-lfdr-simple} as the marginal lfdr since it scores the $i$th null hypothesis as a function of only its $p$-value $p_i$. In practice, we would need to estimate the quantities $\bar{\pi}_0,\bar{f}_0,\bar{f}$ appearing in \eqref{eq:def-lfdr-simple}, so our decision to reject or accept the $i$th null hypothesis eventually depends on all of $p_1,\dots,p_m$. In the absence of further contextual information, it is natural to require the decision rule to be symmetric with respect to the order in which the $p$-values are observed. This symmetry elicits another oracle function, called the compound lfdr, which plays a role parallel to that of the $\mlfdr$ in characterizing the best permutation equivariant decision rule.

We say that a decision rule $\d(\p) \coloneqq (\d_1(\p),\dots,\d_m(\p))$ is is permutation equivariant (\PI) if 
\begin{align}
\label{PI-decision-rule}
    \d(\p)_\pi = \d(\p_\pi) \hspace{1em} \text{ for any } \pi \in \mathcal{S}_m,
\end{align}
where $\mathcal{S}_m$ is the set of permutations on $[m]$, and $v_\pi \coloneqq (v_{\pi(1)},\dots,v_{\pi(m)})$ denotes the vector $v \in \R^m$ permuted by $\pi$. Any multiple testing procedure that uses a rejection threshold which is a function of the order statistics is \PI. For example, the Benjamini-Hochberg procedure \citep{benjamini1995controlling} applied to a list of $p$-values defines a $\PI$ decision rule.

\begin{figure}[t]
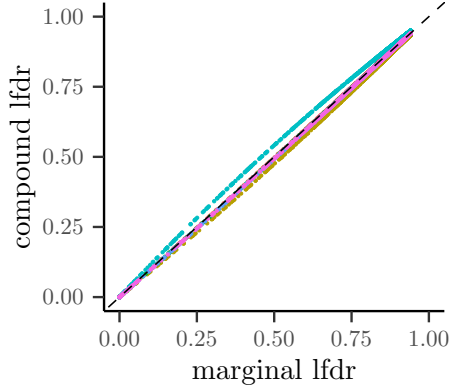

    \centering
    \include{figures/7}
    \vspace*{-30pt}
    \caption{For each of 6 realizations of the vector $(p_1,\dots,p_m)$, with $m=1000$, $\bar{\pi}_0 = 0.8$, $f_0 = 1_{[0,1]}$ and $f_1 = \textnormal{Beta}(1/4,1)$, 
    $\clfdr(\p)$ is approximated numerically and the points $(\mlfdr(p_i),\clfdr_i(\p))$ are plotted with the diagonal $y=x$ shown as a dashed line. Each color represents a different realization of the one-thousand $p$-values.}
    \label{fig:clfdr-vs-simple-lfdr}
\end{figure}

Random shuffling induces an exchangeable Bayesian model:
\begin{equation}
\begin{aligned}
\label{PI-oracle-two-groups-model}
    \pi &\sim \text{Uniform}(\mathcal{S}_m) \\
    \widetilde{H} &\coloneqq H_\pi, \hspace{1em} \widetilde{\p} \coloneqq \p_\pi.
\end{aligned}
\end{equation}
The weighted classification risk of any $\PI$ decision rule $\d$ in this model coincides with its frequentist compound risk, yielding another instance of the fundamental theorem for compound decisions (\cite{zhang2003compound}, \cite{robbins1951asymptotically}, \cite{weinstein2021permutation})
\begin{align*}
    \E L_\lambda(H,\d(\p)) = \widetilde{\E }L_\lambda(\widetilde{H},\d(\widetilde{\p})),
\end{align*}
where $\widetilde{\E}$ marginalizes over $\widetilde{H}$ and $\widetilde{\p}$ generated by \eqref{PI-oracle-two-groups-model}. The right hand side is minimized by the Bayes rule within the exchangeable oracle model \eqref{PI-oracle-two-groups-model}, characterized by the compound lfdr ($\clfdr$),
\begin{align}
\nonumber
    \clfdr_i(t_1,\dots,t_m) &\coloneqq \P(\widetilde{H}_i\text{ is true} \mid \widetilde{\p}=\textbf{t})\\
    \label{eq:clfdr}
    &=\frac{\sum_{\pi \in \mathcal{S}_m: H_{\pi(i)}\text{ is true}} \prod_{j=1}^m f^{(\pi(j))}(t_{j})}{\sum_{\pi \in \mathcal{S}_m} \prod_{j=1}^m f^{(\pi(j))}(t_{j})}, 
\end{align}
for $i=1,\dots,m$ and $\textbf{t} \coloneqq (t_1,\dots,t_m) \in [0,1]^m$. 
It follows that the best $\PI$ decision rule is
\begin{align}
\label{PI-lfdr}
    \d^{*}_i(\p) &\coloneqq \begin{cases}
    1 \hspace{1em}&\text{if } \clfdr_i(\p) \leq \frac{1}{1+\lambda} \\
    0 &\text{else}.
        \end{cases}
\end{align}
This claim follows from a more general relationship between the best $\PE$ decision rule and the Bayes rule with respect to a Haar measure prior (see \cite{eaton2021charles} for a paraphrasing of this result). We also include an elementary proof in appendix \ref{sec-proofs} for completeness.

\begin{figure}[t]
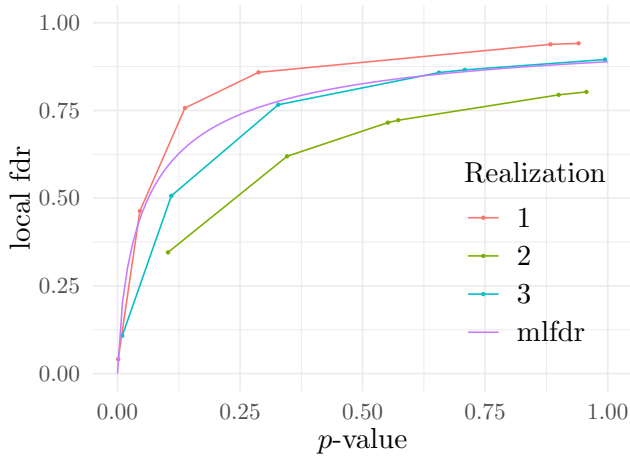

    \centering
    \include{figures/8}
    \vspace*{-30pt}
    \caption{Three realizations of the order statistics of $m=6$ many $p$-values are plotted against their scores $\clfdr_i(\p)$ for $i=1,\dots,6$. In this simulation, $m_0=4$, $f_0 = \text{Uniform}(0,1)$, and $f_1=\text{Beta}(1/4,1)$. The clfdr scores are computed using the realized values $p_{(1)}\leq~\dots \leq~p_{(6)}$ and formula \eqref{eq:clfdr}. The dependence between clfdr scores in any given realization requires that they always sum to $m_0$.} 
    \label{fig:clfdr-interpretation}
\end{figure}

The marginal lfdr is recovered in the exchangeable model \eqref{PI-oracle-two-groups-model} by conditioning on one $p$-value, 
\begin{align*}
    \mlfdr(t) = \P(\widetilde{H}_i \text{ is true} \mid \widetilde{p}_i=t), \hspace{1em} t\in [0,1].
\end{align*} 
Given the true $p$-value densities $f^{(1)},\dots,f^{(m)}$, the $\clfdr$ can typically only be computed in small problems (e.g. $m\leq20$), but can be approximated numerically in larger problems (e.g. $m\approx 1000$) using a method developed by \cite{mccullagh2014asymptotic} for approximating a matrix permanent. Whereas the $\mlfdr$ is a fixed function on $[0,1]$, $\clfdr$ depends on the particular realization of $p$-values, as illustrated in Figure \ref{fig:clfdr-interpretation}. The $\clfdr$ and $\mlfdr$ scores are plotted for six realizations of $p$-values in Figure \ref{fig:clfdr-vs-simple-lfdr}, where they can be seen to roughly coincide for large $m$.

In the next section, we discuss the marginal and compound lfdr functions from a Bayesian perspective. 
Bayesians with an exchangeable prior implicitly report their estimate of the $\clfdr$ via their posterior null probability given all the observations. In light of the previous discussion, this implies that the marginal $\lfdr$ is close to the ``right'' answer in any Bayesian model where the prior is exchangeable and the observations are independent given the truth status of each null hypothesis. 

\subsection{Bayesian interpretation of clfdr and mlfdr}
\label{sec-bayesian-interpretation}

In a Bayesian model, $\clfdr_i(\p)$ is the conditional probability that the $i$th null hypothesis is true, given the data and the empirical distribution of the underlying parameters. For context, suppose there is an exchangeable sequence of latent variables $\theta_1,\ldots,\theta_m$ taking values in some parameter space $\Theta$, and conditional on $\theta=(\theta_1,\ldots,\theta_m)$, the data is drawn according to
\begin{align}
\label{eq:likelihood}
    p_i \mid \theta \sim f_{\theta_i}, \hspace{1em} \text{independently for }i=1,\dots,m.
\end{align}
A standard example is the normal location model, where $f_{\theta_i}(\bar{\Phi}^{-1}(p_i)) = \phi(\bar{\Phi}^{-1}(p_i)-\theta_i)$ and $\phi$ is the standard normal density. The more general setting is recovered by taking the parameter to be $\theta_i = (H_i,f^{(i)})$ and the parameter space to be $\Theta = \{0,1\} \times \{\text{all densities on }[0,1]\}$. 

For a given realization of $\theta$, the marginal and compound lfdr in a Bayesian model with an exchangeable prior on $\theta$ are:
\begin{align}
\label{marginal-lfdr}
    \mlfdr(t;G_m) &= \P(\theta_i=0 \mid p_i=t, G_m), \\
    \label{compound-lfdr}
    \clfdr_{i}(\textbf{t};G_m) &= \P(\theta_{i} =0 \mid \p=\textbf{t},G_m)
\end{align}
where $t \in [0,1]$, $\textbf{t} \in [0,1]^m$, 
\begin{align*}
    G_m(t) \coloneqq m^{-1}\sum_{i=1}^m 1 \{\theta_i \leq t\}
\end{align*} 
is the empirical cumulative distribution function of the true effects, and $\{\theta_i=0\}$ is the null event\footnote{previously denoted ``$H_i\text{ is true}$''}. This definition appears ambiguous, because up until this point, the marginal and compound lfdr have only been defined in a strictly frequentist model. To clarify, conditioned on a specific realization of $\theta$, the joint distribution \eqref{eq:likelihood} defines a frequentist model, and within this model the frequentist $\mlfdr$ and $\clfdr$ functions are equivalent to posterior probabilities that condition also on $G_m$ within the ambient Bayesian model. 
\begin{proposition}
\label{prop-exchangeable-Bayes-interpretation}
    Suppose $\theta_1,\dots,\theta_m$ is an exchangeable sequence of latent variables, and that given $\theta=(\theta_1,\dots,\theta_m)$, the $p$-values are drawn according to \eqref{eq:likelihood}. Let $\mlfdr(t;G_m)$ and $\clfdr_i(\textbf{t};G_m)$ be defined as in \eqref{marginal-lfdr} and \eqref{compound-lfdr} where $\textbf{t}\coloneqq (t_1,\dots,t_m) \in [0,1]^m$. Then 
    \begin{align*}
        \mlfdr(t;G_m) &= \frac{\bar{\pi}_0 f_0(t)}{\frac{1}{m}\sum_{i=1}^m f_{\theta_i}(t)} \\
        \clfdr_i(\textbf{t};G_m) &= \frac{\sum_{\pi \in \mathcal{S}_m:\theta_{\pi(i)}=0} \prod_{j=1}^m f_{\theta_{\pi(j)}}(t_j)}{\sum_{\pi \in \mathcal{S}_m} \prod_{j=1}^m f_{\theta_{\pi(j)}}(t_j)},
    \end{align*}
    for $i=1,\dots,m$, where $\bar{\pi}_0\coloneqq \frac{\# \{i:\theta_i=0\}}{m}$.
\end{proposition}
\noindent For the compound lfdr, there is a large class of Bayesians (essentially, ones with exchangeable priors over $(\theta_i)_{i=1}^m$) for whom their posterior credence in each null hypothesis coincides with their Bayes estimate of compound lfdr. In this sense, we might say Bayesians with exchangeable priors are all in agreement that the compound lfdr is the right quantity to estimate.
The same can nearly be said about the marginal lfdr, for the smaller subclass of Bayesians who look marginally at the data for each hypothesis. For these Bayesians, the posterior probability given a single $p_i$ coincides with their conditional expectation of the $\mlfdr$. These claims are formalized in the next proposition, which is a straightforward consequence of the tower property of conditional expectations.

\begin{proposition}
\label{prop-bayes-estimate-flfdr}

    Suppose the sequence $\{(\theta_i,p_i)\}_{i=1}^m$ is exchangeable and \eqref{eq:likelihood} holds for each $i=1,\dots,m$. Then
    \begin{align}
    \label{eq:exchangeable-prior}
        \P(\theta_i=0 \mid \p) &= \E \left[ \clfdr_i(\p ; \theta) \mid \p \right].
    \end{align}
    Marginally, we have for each $i=1,\dots,m$
    \begin{align}
    \label{eq:cond-iid-prior}
        \P(\theta_i=0 \mid p_i) &= \E \left[ \mlfdr(p_i; \theta) \mid p_i\right].
    \end{align}
    
\end{proposition}

If we can obtain a good estimator of the compound lfdr given structural assumptions like monotonicity, then any Bayesian with an exchangeable prior on the hypotheses should be fairly satisfied with using it to make predictions, since the predictions they would make are just their estimate of the same quantity. In particular, in many large problems, most of these Bayesian observers would converge on similar estimates for compound lfdr. In such cases, a good frequentist estimator of compound lfdr should also give about the same answer.

The marginal lfdr is computationally simpler to evaluate than the compound lfdr, and under sufficiently regular conditions, their ratio tends to 1 as $m \to \infty$. 

\begin{lemma}
\label{prop-asymptotic-PI-rule}
    Suppose $p_i \sim f^{(i)}$ are drawn independently for $i=1,\dots,m$ where each $f^{(i)}$ is a continuous density.
    $f^{(i)}=f_0$ when $H_i=0$ and $f^{(i)}=f_1$ when $H_i=1$. If $\frac{m_0}{m} \to \pi_0 \in (0,1)$ as $m\to\infty$, and $\var\left(\frac{f_1}{f_0}(p_1)\right) \vee \var\left(\frac{f_0}{f_1}(p_2)\right) < \infty$ when $p_1 \sim f_0$ and $p_2 \sim f_1$, then we have
    for each $i=1,2,\dots$
    \begin{align*}
        \P\left(\left|\frac{\clfdr_i(\p)}{\mlfdr(p_i)} - 1 \right| > m^{-1/2} (\log m)^{3/2} \right) \leq \frac{C}{(\log m)^3}
    \end{align*}
    for some constant $C>0$ when $m$ is sufficiently large.
\end{lemma}

\subsection{Discrete-uniform null distribution}
\label{subsec-discrete-uniform}

\begin{figure}[t]
    \centering
    \includegraphics[width=11cm,height=8cm]{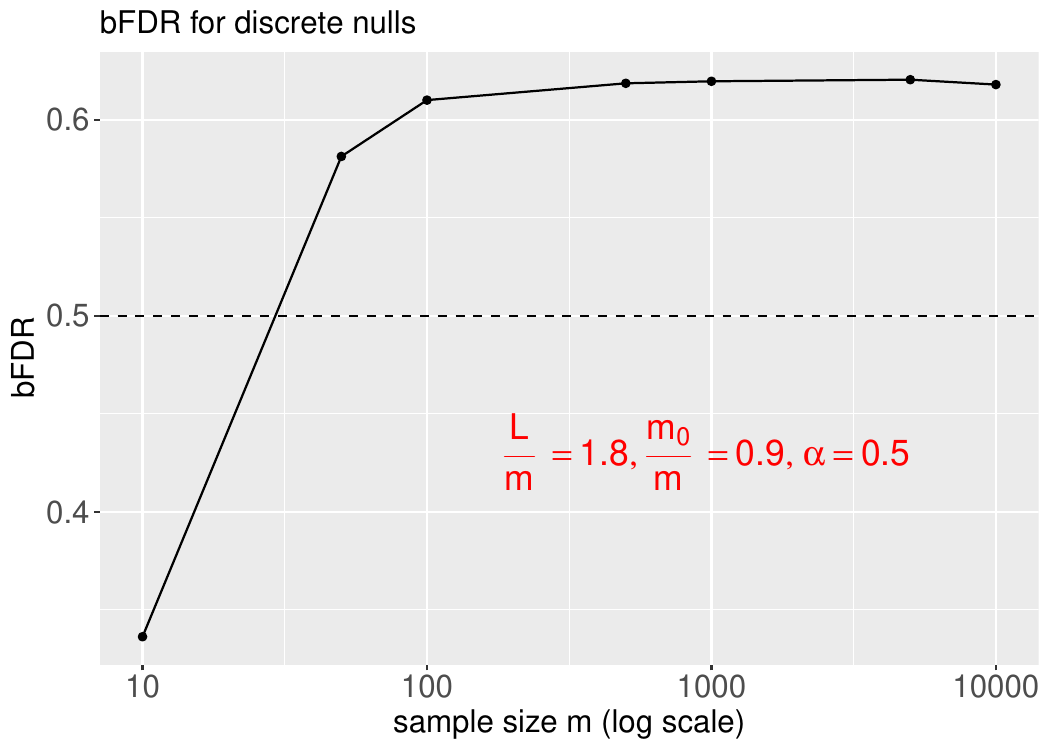}
    \caption{The simulation setting is $\bar{\pi}_0 = 0.9$, $\q=0.5$, and $L=1.8m$, and the non-nulls are fixed along the grid $\left\{\frac{1}{L},\frac{2}{L},\dots,\frac{m_1}{L}\right\}$. The sample size ranges from $m \in \{10,50,100,500,1000,5000,10^4\}$. The bFDR of the SL procedure is estimated using $N=10^5$ Monte Carlo samples.}
    \label{fig:discrete-null-asymp}
\end{figure}

The results that follow in this section are asymptotic, and regard the boundary FDR of the SL procedure run directly on the discrete $p$-values. If either $L \to \infty$ with $m$ fixed, or $m \to \infty$ with $L$ fixed, then the boundary FDR is controlled asymptotically (Theorems \ref{thm-discrete-null} and \ref{thm-discrete-L-fixed}). Things are less clear when $m/L$ converges to a constant. In this case, numerical evidence suggests the bFDR guarantee may be violated even as $m,L$ get large, with $L \approx 1.8m$ (Figure \ref{fig:discrete-null-asymp}). Lemma \ref{lem-discrete-uniform} sheds some light on the interplay between $m$ and $L$ and is important for proving our next result.

\begin{thm}
\label{thm-discrete-null}
    Let $p_1,\dots,p_m$ be independent random variables with the same support, where $p_i \sim \textnormal{Uniform}\left\{\frac{1}{L},\frac{2}{L}\dots,\frac{L-1}{L},1\right\}$ when $H_i$ is true. For any fixed $\q \in [0,1]$ and $m$, we have
    \begin{align*}
        \bFDR(\mathcal{R}_\q) \leq \bar{\pi}_0 \q + O\left(\frac{m^2}{L}\right),
    \end{align*}
    as $L \to \infty$.
\end{thm}
The proof of Theorem \ref{thm-discrete-null} is similar to that of Theorem \ref{main-RFDR-control}, both of which start by splitting the boundary FDR into a symmetric contribution from each null,
\begin{align*}
\bFDR(\mathcal{R}_\q) = m_0 \P(\textnormal{rank}(p_m) = R),
\end{align*}
breaking ties uniformly at random, and assuming without loss of generality that $H_m=0$. As $L \to \infty$, the probability that $p_m$ is the realized boundary rejection is asymptotically no larger than $\q/m$.
\begin{lemma}
\label{lem-discrete-uniform}
    Let $p_1,\dots,p_{m-1} \in [0,1]$ be deterministic (non-random) variables, and suppose that $p_m \sim \textnormal{Uniform}\left\{ \frac{1}{L}, \frac{2}{L},\dots,\frac{L-1}{L},1\right\}$. Then as $L \to \infty$,
    \begin{align*}
        \P(\textnormal{rank}(p_m) = R) \leq \frac{\q}{m} + O\left( \frac{m}{L}\right),
    \end{align*}
    as $L \to \infty$, where $\textnormal{rank}(p_m) \coloneqq \# \{i : p_i \leq p_m\}$ is the rank of $p_m$ among the full list $p_1,\dots,p_m$, breaking ties at random.
\end{lemma}
The assumption that $L \to \infty$ in the previous lemma is needed to rule out small cases where it is possible for the inequality to be violated\footnote{A counterexample is obtained by setting $p_1,\dots,p_{m-1}$ to specific values for which $\P(p_{(R)}=p_m) > \frac{2\q}{m}$ when $L=9,m=6,\q=1/2$. More details about this counter-example can be found in Section \ref{subsubsec-CE-discrete-uniform}.} by a factor of 2. Our next result characterizes the boundary FDR of the SL procedure as $m \to \infty$, keeping $L$ fixed.

\begin{thm}
\label{thm-discrete-L-fixed}
    Let $p_1,\dots,p_m$ be independent random variables with the same support, where $p_i \sim \textnormal{Uniform}\left\{ \frac{1}{L}, \frac{2}{L},\dots,\frac{L-1}{L},1\right\}$ when $H_i$ is true, for some fixed $L$. Suppose as $m \to \infty$ that the average probability mass function $\bar{f}$ converges to a limiting pmf $f^*$, supported on $\left\{ \frac{1}{L}, \frac{2}{L},\dots,\frac{L-1}{L},1\right\}$, and that $\frac{m_0}{m} \to \pi_0^* \in (0,1)$. Further assume that there is a unique maximizer $\ell^*$ of the population-level objective,
    \begin{align*}
        \ell^* \coloneqq \underset{\ell=0,\dots,L}{\operatorname{argmax}} \hspace{.4em}\left\{ \q  \sum_{k=0}^\ell f^*(k/L) -\ell/L\right\}.
    \end{align*}
    Then we have
    \begin{align*}
        \lim_{m \to \infty} \bFDR(\mathcal{R}_\q) = \frac{\pi^*_0}{L f^*(\ell^*/L)}\cdot 1_{\{\ell^*>0\}}.
    \end{align*}
\end{thm}

\begin{remark}
    It follows as a corollary of Theorem \ref{thm-discrete-L-fixed} that $\bFDR(\mathcal{R}_\q) \leq \pi_0^* \q$ asymptotically (as $m \to \infty$, keeping $L$ fixed) because $\ell^*$ must occur at an $\ell \leq L$ for which the discrete difference sequence is non-negative,
\begin{align*}
    \q f^*(\ell^*/L) - 1/L \geq 0 ,
\end{align*}
which implies that $\frac{\pi_0^* /L}{f^*(\ell^*/L)} \leq \pi_0^* \q$.
\end{remark}

\section{Proofs of technical results}
\label{sec-proofs}

\paragraph{Notation.} In what follows, the notation $a_m \sim b_m$ for two sequences $(a_m)$ and $(b_m)$ means that $a_m/b_m \to 1$ as $m \to \infty$. The notation $a_m \asymp b_m$ means there exist constants $c,C>0$ and $M \in \mathbb{N}$ such that for all $m > M$,
\begin{align*}
    c b_m \leq a_m \leq C b_m.
\end{align*}

\begin{proof}[Proof of Theorem \ref{thm:conditional-probability}]
    Given that some $z_J=t$, the index $J$ is a random variable satisfying
    \begin{align*}
        \P(J=j \mid z_J=t) \propto f^{(j)}(t).
    \end{align*}
    By the continuous density assumption, $J$ is almost surely unique. Therefore,
    \[
        \P(J\in\H_0 \mid z_J=t) = \frac{\sum_{j \in \H_0} f^{(j)}(t)}{\sum_{j=1}^m f^{(j)}(t)} = \frac{\bar{\pi}_0 f_0(t)}{\bar{f}(t)}. \qedhere
    \]
\end{proof}

\begin{proof}[Proof of Theorem \ref{thm:calibration}]

    The continuous density assumption implies
    \begin{align*}
        \P\big(\bigcup_{j\in\H_0} \{|\ell_j-\q|\leq \e\}\big) &\sim \sum_{j\in\H_0} \int_{\{t:|\lfdr(t)-\q|\leq \e\}} f_0(t)dt \\
        \P\big(\bigcup_{j=1}^m \{|\ell_j-\q|\leq \e\}\big) &\sim m \int_{\{t:|\lfdr(t)-\q|\leq \e\}} \bar{f}(t)dt,
    \end{align*}
    as $\e \to 0$. Now since $\bar{\pi}_0f_0(t) \sim \q \bar{f}(t)$ for any $t$ such that $|\lfdr(t)-\q| \leq \e$, the ratio tends to $\q$.
    
    Since $g$ is calibrated,
    \begin{align*}
        g(t) &= \E\left[1\{H_I \text{ is true}\} \mid g(z_I)= g(t)\right] \\
        &= \E \left[\mlfdr(z_I) \mid g(z_I) = g(t) \right],
    \end{align*}
    by the tower property, where $I \sim \text{Uniform}\{1,\dots,m\}$.
\end{proof}

\begin{proof}[Proof of Theorem \ref{thm:compound-decisions}]
    The expected weighted classification loss can be re-expressed
    \begin{align*}
        \E L_\lambda(H,\d) = \E_{z_I\sim \bar{f}}  \; \ell_\lambda(H_I,g(z_I)),
    \end{align*}
    where $I\sim \text{Uniform}\{1,\dots,m\}$ and $\ell_\lambda(h,y)$ is the per-instance loss:
    \begin{align*}
        \ell_\lambda(h,y) = \lambda \; 1\{h \text{ is true},y=\text{reject}\} 
        + 1\{h \text{ is false},y=\text{accept}\}.
    \end{align*}
    Since $z_I \sim \bar{f} = \bar{\pi}_0 f_0 + (1-\bar{\pi}_0) \bar{f}_1$ follows a Bayesian two-groups model, the expected loss is minimized by the Bayes rule \citep{sun2007oracle}, which is characterized by the local fdr in this two-groups model,
    \begin{align}
    \label{eq:lfdr-uniform-index}
        \lfdr(t) = \P(H_I\text{ is true} \mid z_I=t).
    \end{align}
\end{proof}

\begin{proof}[Proof of Theorem \ref{pfdr-thm}]
Let $V \equiv V(z-\e,z+\e)$ and $R\equiv R(z-\e,z+~\e)$. First, $\mFDR([t-\e,t+\e])$ is equal to
\begin{align*}
     &= \frac{\E\left[\sum_{i=1}^m (1-H_i) 1\{z_i \in [t-\e,t+\e])\}\right]}{\E\left[ \sum_{i=1}^m 1\{z_i \in [t-\e,t+\e]\}\right]} \\
    &= \frac{\sum_{i=1}^m (1-H_i) \P(z_i \in [t-\e,t+\e])}{\sum_{i=1}^m \P(z_i \in [t-\e,t+\e])} \\
    &\sim \frac{\sum_{i:H_i=0} f^{(i)}(t)\e}{\sum_{i=1}^m f^{(i)}(t) \e} = \mlfdr(t).
\end{align*}
The third line holds because each density $f^{(i)}$ is continuous, so the probability that $z_i \in [t-\e,t+\e]$ is proportional to the density evaluated at some point in this interval, multiplied by the length of the interval.

Next, the pFDR is
\begin{align*}
    \E(V/R \mid R>0) = \frac{\E(V/R \cdot 1\{R>0\})}{\P(R>0)}.
\end{align*}
Since $\P(R=1 \mid R>0) \to 1$ as $\e \to 0$, we have
\begin{align*}
    \E(V/R \cdot 1\{R>0\}) \sim \P\left( \cup_{j \in \H_0} \{|Z_j-z|<\e\}\right).
\end{align*}
Since $\P(R>0) = \P\left(\cup_{j \in [m]} \{|Z_j-z|<\e\}\right)$, the ratio is
\begin{align*}
    &= \lim_{\e \to 0} \; \frac{\P\left(\cup_{j \in \H_0} \{|Z_j-z|<\e\} \right)}{\P\left(\cup_{j\in[m]} \{|Z_j-z|<\e\} \right)} \\
    &= \lim_{\e \to 0} \; \frac{\sum_{j \in \H_0} \P\{|Z_j-z|<\e\} + O(\e^2) }{\sum_{j \in [m]} \P\{|Z_j-z|<\e\} + O(\e^2) } = \lfdr(z).
\end{align*}
\end{proof}

\begin{proof}[Proof of Theorem \ref{thm-calibration-lfdp}]
    Define the count variables
    \begin{align*}
        N_0 &\coloneqq \# \{i \leq m : p_i \in [t,t+\e_m], H_i=0\} \\
        N &\coloneqq \# \{i \leq m : p_i \in [t,t+\e_m] \},
    \end{align*}
    so we may write $\FDP\left(\{H_i : p_i \in [t,t+\e_m], i\leq m\} \right) = N_0/N$. Letting $F_0$ and $\bar{F}_m$ denote the cdfs of $f_0$ and $\bar{f}_m$ respectively, we have 
    \begin{align*}
        \frac{\E N_0}{\E N} &= \frac{m_0 (F_0(t+\e_m)-F_0(t))}{m(\bar{F}_m(t+\e_m)-\bar{F}_m(t))} \sim \frac{\bar{\pi}_{0,m} f_0(t)}{\bar{f}_m(t)} = \lfdr_m(t).
    \end{align*}
    To show \eqref{thm-calib-part1}, it suffices to show
    \begin{align*}
        \frac{N_0}{\E N_0} \stackrel{p}{\to} 1, \hspace{1em}\frac{N}{\E N} \stackrel{p}{\to} 1,
    \end{align*}
    from which \eqref{thm-calib-part1} follows by Slutsky's theorem. By Chebyshev's inequality, for any $\delta > 0$,
    \begin{align*}
        \P\left(\left| \frac{N_0}{\E N_0}-1 \right| > \d \right) \leq \d^{-2} \frac{\var(N_0)}{(\E N_0)^2} \leq \d^{-2}\frac{m_0 (F_0(t+\e_m)-F_0(t))}{m_0^2 (F_0(t+\e_m)-F_0(t))^2} \sim \d^{-2}\frac{1}{m_0 f_0(t) \e_m} \to 0,
    \end{align*}
    for any $t$ with $f_0(t)>0$. By a similar argument,
    \begin{align*}
        \P\left(\left| \frac{N}{\E N}-1 \right| > \d \right) \leq \d^{-2} \frac{1}{m \bar{f}_m(t) \e_m} \leq \d^{-2} \frac{1}{m_0 f_0(t) \e_m} \to 0,
    \end{align*}
    completing the proof of \eqref{thm-calib-part1}. For \eqref{thm-calib-part2}, let $M_0$ and $M$ be defined:
    \begin{align*}
        M_0 &\coloneqq \# \{i \leq m : \lfdr_m(z_i) \in [\alpha,\alpha+\e_m] , H_i=0\} \\
        M &\coloneqq \# \{i \leq m : \lfdr_m(z_i) \in [\alpha,\alpha+\e_m]\}
    \end{align*}
    so that $\FDP\left(\{H_i : \lfdr(z_i) \in [\alpha,\alpha+\e_m], i\leq m\} \right) = M_0/M$. Letting $Z \sim f_0$, we have
    \begin{align*}
    \frac{\E M_0}{\E M} &= \frac{m_0 \P(\lfdr_m(Z) \in [\alpha,\alpha+\e_m])}{\sum_{i=1}^m \P(\lfdr_m(z_i) \in [\alpha,\alpha+\e_m])} \\
    &= \frac{m_0 \int_{\{t:\lfdr_m(t) \in [\alpha,\alpha+\e_m]\}} f_0(t) \; \textnormal{d}t}{m \int_{ \{t: \lfdr_m(t) \in [\alpha,\alpha+\e_m]\}}\bar{f}_m(t)\; \textnormal{d} t},
    \end{align*}
    which implies
    \begin{align*}
        \alpha \leq \frac{\E M_0}{\E M} \leq \alpha + \e_m,
    \end{align*}
    since $m_0 f_0(t) /m = \bar{\pi}_{0,m} f_0(t) \leq \bar{f}(t) (\alpha + \e_m)$ and $m_0 f_0(t) / m \geq \alpha \bar{f}_m(t)$ on the event where $\lfdr_m(t) \in [\alpha,\alpha+\e_m]$. Now it suffices to show
    \begin{align*}
        \frac{M_0}{\E M_0} \stackrel{p}{\to} 1, \hspace{1em} \frac{M}{\E M} \stackrel{p}{\to} 1,
    \end{align*}
    from which \eqref{thm-calib-part2} follows by Slutsky's theorem. By Chebyshev,
    \begin{align*}
        \P \left(\left| \frac{M_0}{\E M_0}-1 \right| > \d \right) \leq \d^{-2}\frac{\var(M_0)}{(\E M_0)^2} \leq \d^{-2} \frac{m_0 \P(\alpha \leq \lfdr_m(Z) \leq \alpha+\e_m)}{m_0^2 \P(\alpha \leq \lfdr_m(Z) \leq \alpha+\e_m)^2},
    \end{align*}
    where $Z \sim f_0$. Suppose $\lfdr_m(Z)$ has density $g_m$. Then
    \begin{align*}
        \P(\alpha \leq \lfdr_m(Z) \leq \alpha+\e_m) = \int_{\alpha}^{\alpha+\e_m} g_m(t) \; \textnormal{d}t \sim g_m(\alpha) \; \e_m,
    \end{align*}
    as $m \to \infty$. Then since $m \e_m \to \infty$ and $g_m(\alpha)$ is bounded away from zero, we have for any $\delta>0$
    \begin{align*}
        \P \left(\left| \frac{M_0}{\E M_0}-1 \right| > \d \right) \leq \d^{-2} \frac{1}{g_m(\alpha) m_0 \e_m} \to 0
    \end{align*}
    as $m \to \infty$. A similar argument shows
    \begin{align*}
        \P \left(\left| \frac{M}{\E M}-1 \right| > \d \right) \leq \d^{-2} \frac{1}{g_m(\alpha) m_0 \e_m} \to 0,
    \end{align*}
    which completes the proof.
\end{proof}

\begin{proof}[Proof of Proposition \ref{prop:uniform-convergence}]
    By Theorem 5.7 in \cite{van2000asymptotic}, it suffices to check that $f^*_m$ is well-separated, i.e. for every $\e > 0$,
    \begin{align*}
        \sup_{f \in \mathcal{F}: D(f^*_m \| f) \geq \e} M_m^*(f) < M_m^*(f^*_m).
    \end{align*}
    For any $f \in \mathcal{F}$ with $D(f^*_m \| f) \geq \e$, we have
    \begin{align*}
        M^*_m(f^*_m) &= \E_{\bar{f}} \log \bar{f}(Z) \\
        &= D(\bar{f} \| f)-D(\bar{f} \| f) + \E_{\bar{f}}\log \bar{f}(Z) \\
        &\geq \e + M_m^*(f),
    \end{align*}
    since $D(\bar{f} \| f) = D(f^*_m\| f) \geq \e$.
\end{proof}

\begin{figure}[t]
    \centering
    \includegraphics[width=.75\textwidth,height=5cm]{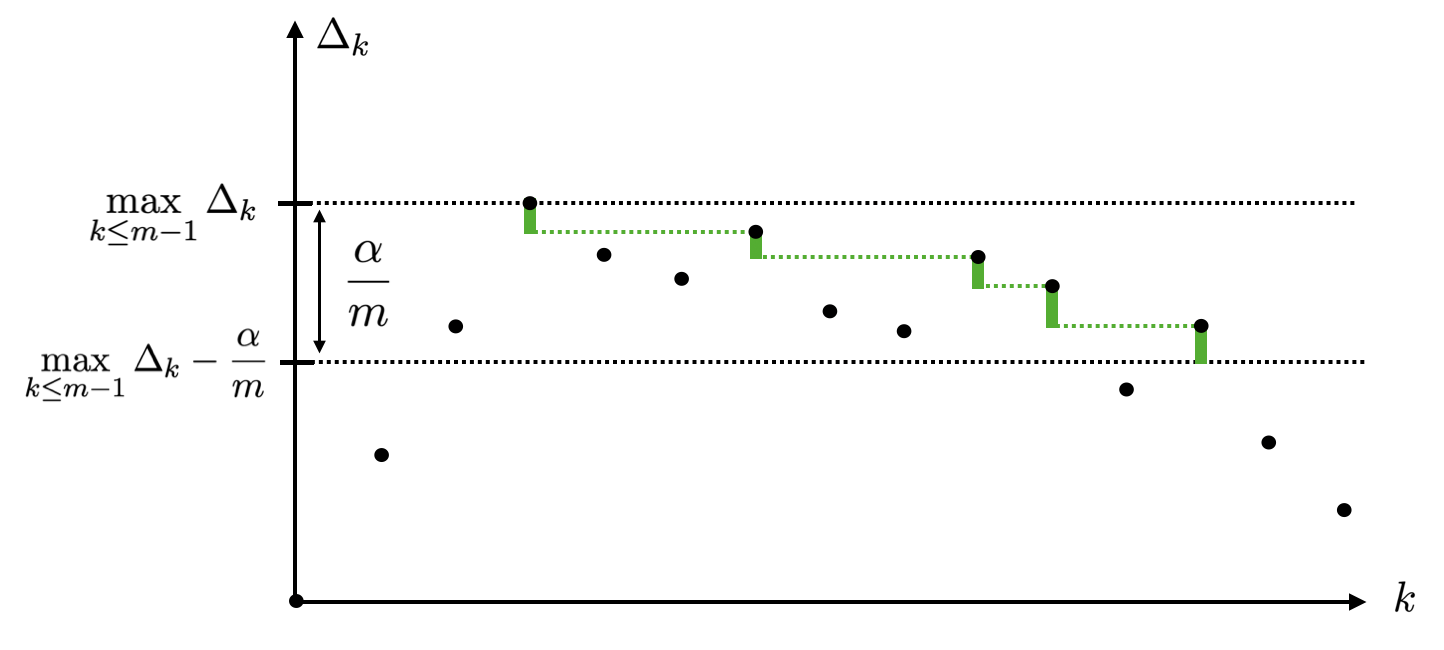}
    \caption{Each length of the interval range in which $p_m$ achieves the maximum in (\ref{num-SL-rejections}) is indicated by a vertical green bar, and the sum of these lengths is $\frac{\q}{m}$.}
    \label{fig:telescope-sum}
\end{figure}

\begin{proof}[Alternative proof of Theorem \ref{main-RFDR-control}]
Suppose without loss of generality that $H_m\text{ is true}$. Then by exchangeability, the $\bFDR$ of the SL method is
\begin{align*}
    \P(H_{(R_\q)}\text{ is true}, R_\q>0) = m\bar{\pi}_0 \P(p_{(R_\q)}=p_m,R_\q>0).
\end{align*} 
Let $q_{(1)} \leq \dots \leq q_{(m-1)}$ denote the order statistics of $p_1,\dots,p_{m-1}$, and note that $p_m$ achieves the maximum in (\ref{num-SL-rejections}) as the $(k+1)$th order statistic if $q_{(k)}< p_m < q_{(k+1)}$ and 
\begin{align*}
    &\frac{\q (k+1)}{m} - p_m > \left[\max_{j=k+1,\dots,m-1} \left\{\Delta_j + \frac{\q}{m}\right\} \right] \vee \left[ \max_{j=0,\dots,k} \Delta_j \right],
\end{align*}
for $k\leq m-1$, where $q_{(0)}\coloneqq 0$ and $\Delta_j \coloneqq \frac{\q j}{m}-~q_{(j)}$. Rearranging the above inequalities gives the range in which $p_m$ achieves the maximum and is equal to the $(k+1)$th order statistic, i.e. $q_{(k)} < p_m$ and
\begin{align*}
    p_m < \frac{\q k}{m} - \left[\max_{j=k+1,\dots,m-1} \Delta_j \right] \vee \left[ \max_{j=0,\dots,k} \Delta_j-\frac{\q}{m} \right].
\end{align*}
This range is non-empty when $\Delta_k$ exceeds each of  $\Delta_{k+1},\dots,\Delta_{m-1}$ as well as $\max_{j=0,\dots,m-1}\Delta_j - \frac{\q}{m}$, and has length 
\begin{align*}
    \Delta_k - \left[\max_{j=k+1,\dots,m-1}\Delta_j\right] \vee \left[\max_{j=0,\dots,k} \Delta_j - \frac{\q}{m}\right]
\end{align*}
The sum of lengths of the non-empty ranges is telescoping and equal to $\frac{\q}{m}$, as illustrated in Figure~\ref{fig:telescope-sum}. 
\end{proof}

\begin{proof}[Proof of Proposition \ref{prop-exp-fam}]
    When $Z \sim g_\theta$, the density of $p = 1-G_{\theta_0}(Z)$ is
    \begin{align*}
        \frac{\de}{\de t} \P_\theta(p \leq t) = \frac{g_\theta}{g_{\theta_0}}(G^{-1}_{\theta_0}(1-t)).
    \end{align*}
    At $\theta=\theta_0$, the above ratio is equal to 1. When $\theta \leq \theta_0$, the log density has a positive derivative in $\theta$ when
    \begin{align*}
        \frac{\de}{\de \theta}\left[\log \frac{g_\theta}{g_{\theta_0}}(G_{\theta_0}^{-1}(1-t))\right] = G_{\theta_0}^{-1}(1-t) - \E_\theta(Z) > 0
    \end{align*}
    which holds for all $t \leq \q$ if $G_{\theta_0}^{-1}(1-\q) > \E_{\theta_0}(Z)$.
\end{proof}
\begin{proposition}
\label{prop-best-PE-rule}
    In the setting of section \ref{subsec-compound-lfdr}, the best $\PE$ decision rule for minimizing the weighted classification risk is defined by \eqref{eq:clfdr} and \eqref{PI-lfdr}.
\end{proposition}
\begin{proof}
    For any $\PE$ decision rule $\d$,
    \begin{align*}
        \widetilde{\E} L_\lambda (\widetilde{H},\d(\widetilde{\p})) &= \frac{1}{m!} \sum_{\sigma \in \mathcal{S}_m} \E L_\lambda(H_\sigma, \d(\p_\sigma)) \\
        &= \frac{1}{m!} \sum_{\sigma \in \mathcal{S}_m} \E L_\lambda(H_\sigma, \d(\p)_\sigma) \\
        &= \E L_\lambda(H,\d(\p)).
    \end{align*}
    The Bayes rule
    \begin{align*}
        \d^* = \underset{\d}{\operatorname{argmin}} \hspace{.4em}  \widetilde{\E} \big[ L_\lambda (\widetilde{H},\d(\widetilde{\p}))  \mid \widetilde{\p} \big],
    \end{align*}
    is itself PE due to exchangeability of $(\widetilde{H}_i,\widetilde{\p}_i)$ across $i=1,\dots,m$. To see this, note that for any $\sigma\in \mathcal{S}_m$, we have 
    \begin{align*}
        (\widetilde{H},\widetilde{\p}) \stackrel{(d)}{=} (\widetilde{H}_\sigma,\widetilde{\p}_\sigma),
    \end{align*} 
    which implies the posterior probability mass function of $\widetilde{H} \mid \widetilde{\p}=t$ at $\widetilde{h} \in \{0,1\}^m$ is equal to the posterior pmf of $\widetilde{H}_\sigma\mid \widetilde{\p}_\sigma=t_\sigma$ at $\widetilde{h}_\sigma$. Thus,
    \begin{align*}
        \d^*(t) &= \underset{h \in \{0,1\}^m}{\operatorname{argmin}} \hspace{.4em}  \widetilde{\E} \big[ L_\lambda (\widetilde{H}_\sigma,h)  \mid \widetilde{\p}_\sigma=t \big] \\
        &= \underset{h \in \{0,1\}^m}{\operatorname{argmin}} \hspace{.4em}  \widetilde{\E} \big[ L_\lambda (\widetilde{H},h_{\sigma^{-1}})  \mid \widetilde{\p}_\sigma=t \big] \\
        &= \bigg[\underset{g \in \{0,1\}^m}{\operatorname{argmin}} \hspace{.4em}  \widetilde{\E} \big[ L_\lambda (\widetilde{H},g)  \mid \widetilde{\p}=t_{\sigma^{-1}} \big] \bigg]_{\sigma}\\
        &= \d^*(t_{\sigma^{-1}})_\sigma.
    \end{align*}
    Since the above holds for any permutation $\sigma$, the Bayes rule in model \eqref{PI-oracle-two-groups-model} is a PE decision rule. Since the average risk in the Bayes model \eqref{PI-oracle-two-groups-model} is equal to the risk function in the frequentist model for every configuration of truth values $H \in \{0,1\}^m$, the Bayes rule is equal to the best $\PE$ rule
    \[
        \d^* = \underset{\d\hspace{.25em}\PE}{\operatorname{argmin}} \hspace{.4em}  \E L_\lambda(H,\d(\p)). \qedhere
    \]
\end{proof}

\begin{proof}[Proof of Proposition \ref{prop-exchangeable-Bayes-interpretation}]
According to Bayes rule, $\P(\theta_i=0 \mid p_i=t, G_m)$ is equal to
\begin{align*}
    &= \frac{\P(\theta_i=0\mid G_m) f_{0}(t)}{\sum_{k=1}^m f_{\theta_{(k)}}(t) \P(\text{rank}(\theta_i)=k \mid G_m)},
\end{align*}
where $\rank(\theta_i)=k$ when $\theta_j<\theta_i$ exactly $k-1$ indices $j \in [m]$, and $\theta_{(1)}\leq \dots \leq \theta_{(m)}$ are the ordered values of $\theta_1,\dots,\theta_m$. Since $\theta_1,\dots,\theta_m$ are exchangeable, the above is equal to
\begin{align*}
    \P(\theta_i=0 \mid p_i=t,G_m) &= \frac{G_m(\{0\}) f_0(t)}{ \frac{1}{m} \sum_{j=1}^m f_{\theta_j}(t)}.
\end{align*}
For \eqref{compound-lfdr}, note that when $p_i = p_{(k)}$, exchangeability implies
\begin{align*}
    \P(&\theta_i=0 \mid p_1=t_1,\dots,p_m=t_m, G_m)\\ &\propto \sum_{\pi \in \mathcal{S}_m:\theta_{\pi(i)}=0} \prod_{j=1}^m f_{\theta_{\pi(j)}}(t_j).
\end{align*}    
\end{proof}

\begin{proof}[Proof of Lemma \ref{prop-asymptotic-PI-rule}]
    The argument is adapted from Theorem 3.1 in \cite{greenshtein2009asymptotic}. Supposing without loss of generality that $H_1=0$ and $H_2=1$,
    \begin{align*}
    \clfdr_i(\p) &= \frac{\bar{\pi}_0 f_0(p_i)}{\bar{\pi}_0 f_0(p_i) + \bar{\pi}_1 f_1(p_i) \cdot X_i},
    \end{align*}
    where $X_i$ is a likelihood ratio,
    \begin{align*}
        X_i \coloneqq \frac{\sum_{\sigma \in \mathcal{S}_m : \sigma(i)=2} \prod_{j \in [m]\backslash \{i\}}^m f^{(\sigma(j))}(p_{j})}{\sum_{\sigma \in \mathcal{S}_m : \sigma(i)=1} \prod_{j \in [m]\backslash \{i\}}^m f^{(\sigma(j))}(p_{j})}
    \end{align*}
    for testing between the following two hypotheses:
    \begin{align*}
        &\text{Hyp}_{0} : \text{Observe a random permutation of $p_{-i}$ when $H_i = 1$} \\
        &\text{Hyp}_{1} : \text{Observe a random permutation of $p_{-i}$ when $H_i = 0$},
    \end{align*}
    where the permutations are drawn uniformly at random from $\mathcal{S}_{m-1}$. A simpler testing problem is:
    \begin{align*}
        &\widetilde{\text{Hyp}}_{0} : \widetilde{p}_1,\dots,\widetilde{p}_{m_0} \stackrel{\iid}{\sim} f_0, \text{ and } (\widetilde{p}_{m_0+1},\dots,\widetilde{p}_{m-1}) \stackrel{\iid}{\sim} f_1 \\ 
        &\widetilde{\text{Hyp}}_{1} : \frac{1}{m_0} \sum_{\ell=1}^{m_0} \left[ (\widetilde{p}_{1:m_0})_{-\ell} \stackrel{\iid}{\sim} f_0, (\widetilde{p}_\ell,\widetilde{p}_{m_0+1},\dots,\widetilde{p}_{m-1}) \stackrel{\iid}{\sim} f_1 \right],
    \end{align*}
    since $\text{Hyp}_{0},\text{Hyp}_{1}$ can be obtained from $\widetilde{\text{Hyp}}_{0},\widetilde{\text{Hyp}}_{0}$ by adding a random permutation. If $H_i=0$ (resp. $H_i=1$), then the distribution of $X_i$ is as if the data were generated by $\text{Hyp}_{1}$ (resp. $\text{Hyp}_{0}$). The likelihood ratio of $\widetilde{\text{Hyp}_{1}}$ to $\widetilde{\text{Hyp}_{0}}$ has variance
    \begin{align*}
        \var_0\left(\frac{1}{m_0} \sum_{\ell =1}^{m_0} \frac{f_1}{f_0}(\widetilde{p}_\ell) \right) = \frac{1}{m_0} \var_0\left( \frac{f_1}{f_0}(p_1) \right) \to 0
    \end{align*}
    by assumption, where $\var_0$ denotes the variance operation when $\widetilde{\text{Hyp}}_0$ holds. It follows from Lemma 2.1 in \cite{greenshtein2009asymptotic} that
    \begin{align*}
        \E_{H_i=1} (X_i-1)^2 \leq \widetilde{\E}_0 \left( \frac{1}{m_0} \sum_{\ell =1}^{m_0} \frac{f_1}{f_0}(\widetilde{p}_\ell)-1\right)^2 \to 0.
    \end{align*}
    A symmetric argument yields
    \begin{align*}
        \E_{H_i=0} (X_i-1)^2 \leq \widetilde{\E}_1 \left( \frac{1}{m_1} \sum_{\ell =m_0}^{m-1} \frac{f_0}{f_1}(\widetilde{p}_\ell)-1\right)^2 \to 0,
    \end{align*}
    under the condition that $\var\left( \frac{f_0}{f_1}(p_2) \right)<\infty$ when $H_2=1$. Here we are abusing notation by writing the index $\ell$ from $m_0$ to $m-1$, to denote summing over the $m_1-1$ many $p$-values drawn from $f_1$ in the scenario described by $\widetilde{\text{Hyp}_1}$. It now follows from Chebyshev's inequality that
    \begin{align*}
        &\P\left( \left|\frac{\clfdr_i(\p)}{\mlfdr(p_i)}-1\right| > m^{-1/2}(\log m)^{3/2} \right) \\
        &\leq \P\left(|X_i-1|>m^{-1/2}(\log m)^{3/2}\right) \\
        &\leq \frac{m}{(\log m)^3} \cdot \frac{\var\left( \frac{f_0}{f_1}(p_2) \right) \vee \var\left( \frac{f_1}{f_0}(p_1) \right)}{m_0 \wedge m_1} \\
        &\leq \frac{C}{(\log m)^3}
    \end{align*}
    for some constant $C>0$ as $m \to \infty$, since $\bar{\pi}_0$ is bounded away from zero and one.
\end{proof}

\begin{proof}[Proof of Lemma \ref{lem-discrete-uniform}.]
    By the law of total probability,
    \begin{align*}
        \P(\textnormal{rank}(p_m) = R) &= \sum_{\ell=1}^L \frac{1}{L} \cdot \P(\textnormal{rank}(p_m) = R \mid p_m = \ell/L) \\
        &= \sum_{\ell=1}^L \frac{1}{L} \cdot 1_{\{\hat{\tau}_\q(\ell/L) = \ell / L\}} \cdot \frac{1}{n_{\ell}+1},
    \end{align*}
    where $n_\ell \coloneqq \# \{i < m : p_i = \ell/L\}$, and we have used explicit notation $\hat{\tau}_\q(p_m)$ to denote the threshold $\hat{\tau}_\q$ as a function of $p_m$, 
    \begin{align*}
        \hat{\tau}_\q(p_m) &\coloneqq \underset{p_{(k)}}{\operatorname{argmax}} \hspace{.4em} \left\{ \frac{\q k}{m} - p_{(k)} \right\} \\
        &= \frac{1}{L}\cdot \underset{\ell = 0,\dots,L}{\operatorname{argmax}} \hspace{.4em} \left\{ \frac{\q L}{m} \cdot \#\{i\leq m:p_i \leq \ell/L\} - \ell \right\},
    \end{align*}
    treating $p_1,\dots,p_{m-1}$ as non-random elements of the grid $\{1/L,\dots,L/L\}$. Define
    \begin{align*}
        \Delta_{\ell} \coloneqq \frac{\q L}{m}N_\ell - \ell, \hspace{1em} \ell = 1,\dots,L
    \end{align*}
    where $N_\ell \coloneqq \# \{i<m : p_i \leq \ell/L\}$, and let $\ell^*\coloneqq \underset{\ell}{\operatorname{argmax}} \hspace{.4em}\Delta_{\ell}$. We claim that
    \begin{align}
    \label{discrete-winners-claim}
        \hat{\tau}_\q (\ell/L) = \ell/L \iff \Delta_\ell > \left[ \Delta_{\ell^*} - \frac{\q L}{m} \right] \vee \max_{k>\ell} \Delta_k,
    \end{align}
    which follows from the same argument as in the alternative proof of Theorem \ref{main-RFDR-control} in Appendix \ref{sec-proofs}.
    The claim implies
    \begin{align*}
        \P(\textnormal{rank}(p_m) = R) &=  \sum_{\ell=1}^L \frac{1}{L} \cdot 1_{\{\hat{\tau}_\q(\ell/L) = \ell / L\}} \cdot \frac{1}{n_{\ell}+1} \\
        &\leq \frac{1}{L} \sum_{\ell=1}^L  1_{\left\{\Delta_\ell > \left[ \Delta_{\ell^*} - \frac{\q L}{m} \right] \vee \max_{k>\ell} \Delta_k \right\}} .
    \end{align*}
    Without loss of generality, suppose that $p_1 \leq p_2 \leq \dots \leq p_{m-1}$. If $p_i<p_j$ are among these $p$-values and both $\Delta_{ p_i L },\Delta_{ p_j L }$ satisfy the rhs of \eqref{discrete-winners-claim}, then as $L\to \infty$, it follows from the definition of $\Delta_\ell$ that as $L$ grows,
    \begin{align*}
        \Delta_{\ell} - \Delta_{\ell+1} = 1 \hspace{1em} \text{for all but a fixed number of } \ell \in (Lp_i,Lp_{j}).
    \end{align*}
    Thus for the distinct values $\ell^* =: \ell_1 < \ell_2 < \dots < \ell_k$ for which the rhs of \eqref{discrete-winners-claim} holds, we must have $\Delta_{\ell_{i}}-\Delta_{\ell_{i+1}} = 1$ for all but a fixed number of winners $\ell_i$, the number of which is less than $m$. Therefore
    \begin{align*}
        \P(\textnormal{rank}(p_m) = R) 
        &\leq \frac{1}{L} \sum_{\ell=1}^L  1_{\left\{\Delta_\ell > \left[ \Delta_{\ell^*} - \frac{\q L}{m} \right] \vee \max_{k>\ell} \Delta_k \right\}} \\
        &\leq \frac{1}{L}\left(m + \sum_{i=1}^{k-1} (\Delta_{i}-\Delta_{i+1}) + \Delta_{k}-\Delta_{\ell^*}+\frac{\q L}{m} \right)\\
        &= \frac{1}{L}\left(m + \frac{\q L}{m} \right) = \frac{\q}{m} + O\left(\frac{m}{L}\right),
    \end{align*}
    as $L \to \infty$.
\end{proof}

\begin{proof}[Proof of Theorem \ref{thm-discrete-L-fixed}]
    By symmetry of the nulls,
    \begin{align*}
        \P(H_{(R_\q)}=0) = m_0 \P(\text{rank}(p_m) = R).
    \end{align*}
    By the argument in Theorem \ref{thm-discrete-null}, the probability that $p_m$ is the $R^{\th}$ smallest $p$-value (breaking ties uniformly at random) is
    \begin{align*}
        \P(\text{rank}(p_m) = R) &= \E\left( \frac{1}{L} \sum_{\ell=1}^L 1_{\left\{\Delta_\ell > \left[\Delta_{\ell^*_m}- \q L/m \right] \vee \max_{j > \ell} \Delta_j \right\}} \cdot \frac{1}{n_{\ell}+1} \right),
    \end{align*}
    where $n_\ell \coloneqq \# \{i < m : p_i = \ell/L\}$ and $\Delta_\ell$ are defined
    \begin{align*}
    \Delta_\ell \coloneqq \frac{\q L}{m} N_\ell - \ell, \hspace{1em} \ell=0,\dots,L,
    \end{align*}
    with $N_\ell \coloneqq \sum_{k=1}^\ell n_k$ and $\ell^*_m \coloneqq \underset{\ell = 0,\dots,L}{\operatorname{argmax}}\hspace{.4em} \Delta_\ell$. As $m \to \infty$, we have the following convergence in probability,
    \begin{align*}
        \frac{n_\ell}{mf^*(\ell/L)} \stackrel{p}{\to} 1, \hspace{1em} \Delta_\ell \stackrel{p}{\to} \q L \sum_{k=0}^\ell f^*(k/L)-\ell.
    \end{align*}
    Since the maximizer of $\q L\sum_{k=0}^\ell f^*(k/L)-\ell$ is unique, we have $\ell^*_m = \ell^*$ with high probability as $m \to \infty$, and
    \begin{align*}
         mf^*(\ell^*/L) \cdot \sum_{\ell=1}^L 1_{\left\{\Delta_\ell > \left[\Delta_{\ell^*_m}- \q L/m \right] \vee \max_{j > \ell} \Delta_j \right\}} \cdot \frac{1}{n_{\ell}+1} \stackrel{\P}{\longrightarrow} 1_{\{\ell^*>0\}},
    \end{align*}
    from which it follows that
    \begin{align*}
        \E \left[\frac{m_0}{mf^*(\ell^*/L)} \cdot \frac{L}{L} \cdot mf^*(\ell^*/L) \cdot \sum_{\ell=1}^L 1_{\left\{\Delta_\ell > \left[\Delta_{\ell^*_m}- \q L/m \right] \vee \max_{j > \ell} \Delta_j \right\}} \cdot \frac{1}{n_{\ell}+1} \right] \to \frac{\pi_0^*}{Lf^*(\ell^*/L)},
    \end{align*}
    when $\ell^*>0$, and zero otherwise.
\end{proof}

\subsubsection{Counterexample to bFDR control under the discrete-uniform assumption}
\label{subsubsec-CE-discrete-uniform}

The SL guarantee breaks down in the setting where each null $p$-value is uniformly distributed on the grid $\left\{\frac{1}{L},\frac{2}{L},\dots,\frac{L-1}{L},1 \right\}$.

\begin{proof}[Counterexample]
    Let $m=6,L=9,\q=1/2$, and the alternative $p$-values are 
    \begin{align*}
        p_1&=p_2 = 1/L\\
        p_i &= i/L \hspace{1em} i=2,3,4.
    \end{align*}
    Then the probability that $p_m \sim \text{Uniform}\left\{\frac{1}{L},\frac{2}{L},\dots,\frac{L-1}{L},1 \right\}$ is the last SL$(\q)$ rejection is
    \begin{align*}
        \P(p_{(R_\q)}=p_m) &= \sum_{\ell=1}^L \frac{1}{L} \P(p_{(R_\q)}=p_m \mid p_m = \ell/L) \\
        &= \sum_{\ell=1}^L \frac{1}{L} \cdot 1_{\{p_{(R_\q)}=\ell/L\}} \cdot \frac{1}{n_\ell + 1},
    \end{align*}
    where $n_{\ell} \coloneqq \# \{i < m : p_i = \ell/L\}$ for $\ell = 1,\dots,L$. It is straightforward to check that $1_{\{p_{(R_\q)}=\ell/L\}}=1$ for $\ell = 1,2,3,4$, so that the above evaluates to
    \begin{align*}
        \P(p_{(R_\q)}=p_m) = \frac{1}{9}\left( \frac{1}{3} + \frac{1}{2} + \frac{1}{2} + \frac{1}{2} \right) =0.204 > 0.167 = \frac{2\q}{m}.
    \end{align*}
\end{proof}

\subsubsection{Counterexample to bFDR control under the super-uniform assumption}
\label{subsubsec-CE-super-uniform}

The SL guarantee may break down when the nulls are super-uniformly distributed instead of Uniform$(0,1)$ distributed, as illustrated by the following example.
\begin{proof}[Counterexample]
    Let $m=2$, $H_1=0$, and $H_2=1$, with $p_1 \sim f^{(1)}$ defined
\begin{align*}
    f^{(1)}(t) \coloneqq \begin{cases}
        \frac{1}{2} \hspace{1em} &0\leq t \leq \frac{1}{4} \\
        \frac{3}{2} &\frac{1}{4} < t \leq \frac{1}{2} \\
        1 &\frac{1}{2} < t \leq 1,
    \end{cases}
\end{align*}
and $p_2 \equiv \frac{1}{4}$. It is straightforward to verify that $f^{(1)}$ is super-uniform, i.e. 
\begin{align*}
    \int_0^t f^{(1)}(s) \de s \leq t 
    \hspace{1em} \text{ for any } t\in [0,1].
\end{align*}
The event $p_{(R_\q)}=p_1$ is equivalent to
\begin{align*}
    \left(\{p_1 \leq p_2\} \cap \left\{p_1 - \frac{\q}{2} < (p_2 - \q) \wedge 0 \right\}\right) \bigcup \left( \{p_1 > p_2\} \cap \left\{ p_1-\q < \left(p_2 - \frac{\q}{2}\right)  \wedge 0\right\} \right).
\end{align*}
Plugging in $\q = \frac{1}{2}$ and $p_2 = \frac{1}{4}$ gives
\begin{align*}
    H_{(R_\q)}=0 \iff p_{(R_\q)}=p_1 \iff p_1 \in \left(1/4,1/2 \right),
\end{align*}
which occurs with probability $\frac{3}{2}\times \frac{1}{4} = \frac{3}{8} > \frac{\q}{2} = \frac{1}{4}$.
\end{proof}

\begin{proof}[Proof of Theorem \ref{flfdr-control}]
Lemma \ref{lemma-threshold-concentration} implies that for $m$ large enough, we have with probability $\geq 1-\d$ that
\begin{align}
\label{whp-thresholds}
\hat{\tau}_{\q}-\tau^*_{\q} \leq \e \coloneqq C' m^{-1/3}\log(2m/\d),
\end{align}
for some constant $C'>0$ depending on $\q,\L$ and $\d$. Since $f$ is decreasing on $(0,\tau^*_\q+\q)$ and has derivative greater than $-\L^{-1}$ over the interval $(\tau^*_{\q},\tau^*_\q+\e)$, the above inequality implies 
\begin{align*}
f(\hat{\tau}_{\q}) \geq f(\tau_{\q}^*) - \L^{-1} \e = \q^{-1} - \L^{-1} \e.
\end{align*}
It follows that
\begin{align*}
    \mlfdr(\hat{\tau}_\q) = \frac{\pi_0}{f(\hat{\tau}_\q)} \leq \pi_0\q + C m^{-1/3}\log(m/\d),
\end{align*}
for another constant $C>0$ depending on $\q,\L$ and $\d$. The other direction follows similarly from Lemma \ref{lemma-threshold-concentration-lower}.
\end{proof}

\begin{lemma}
\label{lemma-switch}
Let $\tau^*_\q$ be a solution to $f(\tau^*_\q)=\q^{-1}$ and let $\hat{\tau}_\q$ denote the rejection threshold of the SL$(\q)$ procedure (\ref{num-SL-rejections}). If $\hat{\tau}_{\q} >\tau^*_{\q} +\e$, then there exists an index $k\geq 1$ for which
\begin{align*}
p_{(i^*+k)} \leq \tau^*_{\q} +\frac{\q k}{m} \hspace{1em}\text{and} \hspace{1em} k > \frac{m\e}{\q},
\end{align*}
where $i^* \coloneqq \max\{i: p_{(i)} \leq \tau^*_{\q}\}$ and $i^*=0$ if no such $i$ exists. 
\end{lemma}

\begin{proof}
Let $\hat{k}$ be the index for which $\hat{\tau}_{\q} =p_{(i^*+\hat{k})}$. The first inequality can be written
\begin{align*}
\frac{i^*+\hat{k}}{m}-\frac{i^*}{m}-\q^{-1}(p_{(i^*+\hat{k})}-\tau^*_{\q}) \geq 0,
\end{align*}
which holds because $F_m(t)-\q^{-1}F_0(t)$ is maximized at $t=p_{(i^*+\hat{k})}$. Since $\hat{\tau}_{\q} >\tau^*_{\q}+\e$, the above inequality implies $\hat{k} > \frac{m\e}{\q}$.
\end{proof}

\begin{lemma}
\label{lemma-threshold-concentration}
Let $\tau^*_\q$ and $\hat{\tau}_\q$ be defined as in Lemma \ref{lemma-switch}, let $\d>0$ and suppose $\bar{f}$ is decreasing on $[\tau^*_\q,\tau^*_\q+\q]$ and that there exists some $\L>0$ for which $\L \leq |\bar{f}'(t)|\leq \L^{-1}$ for all $t$ with $|t-\tau^*_\q| \leq \e$, where $\e \coloneqq \left(\frac{24}{\q L^2}\right)^{1/3} m^{-1/3}\log(2m/\d)$. Then
\begin{align*}
\P(\hat{\tau}_\q > \tau^*_\q+\e) \leq \d, 
\end{align*}
for any $m \geq C(\q,\L,\d)$, a constant depending only on $\q,\L$ and $\d$.
\end{lemma}

\begin{proof}[Proof of Lemma \ref{lemma-threshold-concentration}]

Applying Lemma \ref{lemma-switch} with $\e$ defined as above, we have
\begin{align}
\nonumber
\P(\hat{\tau}_\q > \tau^*_\q+\e) &\leq \sum_{k>\frac{m\e}{\q}} \P\left(p_{(i^*+k)} \leq \tau^*_\q + \frac{\q k}{m} \right) \\
\label{switch-union-bound}
&= \sum_{\frac{m\e}{\q}<k\leq \frac{m\e \log m}{\q}} \P\left( N_k \geq k \right)+\sum_{k> \frac{m\e \log m}{\q}} \P\left( N_k \geq k \right),
\end{align}
where $i^*$ is defined in Lemma \ref{lemma-switch}, and $N_k$ is the number of $p$-values between $\tau^*_\q$ and $\tau^*_\q+~\frac{\q k}{m}$, distributed Generalized-Binomial with sample size $m$ and average success probability $\bar{F}(\tau^*_\q+~\q k/m)-\bar{F}(\tau^*_\q)$,
\begin{align*}
N_k = \sum_{j=1}^m 1_{\{p_j \in (\tau^*_\q,\tau^*_\q +\q k/m)\}} \Rightarrow \E N_k = m(\bar{F}(\tau^*_\q+\q k/m)-\bar{F}(\tau^*_\q)),
\end{align*}
where $\bar{F} \coloneqq \frac{1}{m}\sum_{i=1}^m F^{(i)}$ is the average cdf of the $p$-values. Note that since $\bar{F}' = \bar{f}$, we have by the mean value theorem that
\begin{align*}
\E N_k = m(\bar{F}(\tau^*_\q+\q k/m)-\bar{F}(\tau^*_\q)) = m \bar{f}(\xi) \cdot \frac{\q k}{m} ,
\end{align*}
for some $\xi\in (\tau^*_\q,\tau^*_\q+\q k/m)$. By the monotonicity assumption, $\bar{f}(\xi) \leq \bar{f}(\tau^*_\q)=\q^{-1}$ implies we have $\E N_k \leq k$. Consider the corresponding Binomial random variable, $\widetilde{N}_k \sim \text{Binomial}(m, \bar{F}(\tau^*_\q+\q k/m)-\bar{F}(\tau^*_\q))$. Since $\E \widetilde{N}_k = \E N_k \leq k$, it follows from Theorem 5 in \cite{hoeffding1956distribution} that 
\begin{align*}
\P\left( N_k \geq k \right) \leq \P\left( \widetilde{N}_k \geq k \right) = \P\left(\widetilde{N}_k \geq \E \widetilde{N}_k \cdot \frac{k}{\E \widetilde{N}_k} \right).
\end{align*}
To bound the probability on the right hand side, we use the following bounds on the expectation $\E \widetilde{N}_k = m (\bar{F}(\tau^*_\q+\q k/m)-\bar{F}(\tau^*_\q))$,
\begin{align}
\label{Nk-upper-bound}
\E \widetilde{N}_k &\leq k + \frac{\L m\e^2}{2} - \L\q k\e \\
\label{Nk-lower-bound}
\E \widetilde{N}_k &\geq  \frac{m\e}{2\q}.
\end{align}
Before proving inequalities (\ref{Nk-upper-bound}) and (\ref{Nk-lower-bound}), we show how they can be used to complete the proof. When $\frac{m\e}{\q}< k \leq \frac{m\e \log m}{\q}$, the upper bound (\ref{Nk-upper-bound}) gives $\E \widetilde{N}_k \leq k + \frac{\L m\e^2}{2}-\L\q \e \cdot \frac{m\e}{\q}$, which implies
\begin{align*}
\P\left(\widetilde{N}_k \geq \E \widetilde{N}_k \cdot \frac{k}{\E\widetilde{N}_k}\right) &\leq \P\left(\widetilde{N}_k \geq \E \widetilde{N}_k \cdot \frac{k}{k-\frac{\L m\e^2}{2}}\right).
\end{align*}
Now since $\frac{1}{1-x} \geq 1+x$, the rhs of the above is
\begin{align*}
&\leq \P\left(\widetilde{N}_k \geq \E \widetilde{N}_k \cdot \left( 1+\frac{\L m \e^2}{2k}\right)\right) \leq \exp\left( -\frac{1}{3} \cdot \E \widetilde{N}_k \cdot  \left( \frac{\L m\e^2}{2k}\right)^2\right),
\end{align*}
where the last inequality follows from a Binomial tail bound, recorded in Lemma \ref{lemma-multiplicative-chernoff}. Now using $k\leq \frac{m\e \log m}{\q}$ and applying the lower bound (\ref{Nk-lower-bound}), we obtain
\begin{align*}
&\leq \exp\left( -\frac{1}{3} \cdot \frac{m\e}{2\q} \cdot \left( \frac{\L\q\e}{2\log m} \right)^2\right).
\end{align*}
Simplifying, we have shown that when $\frac{m\e}{\q}<k\leq \frac{m\e \log m}{\q}$, 
\begin{align*}
\P\left( \widetilde{N}_k \geq k\right) \leq \exp\left( -\frac{\q \L^2 m\e^3}{24\log^2 m} \right).
\end{align*}
Plugging in the formula for $\e$, the above inequality implies that the first piece of (\ref{switch-union-bound}) is bounded,
\begin{align}
\label{small-k-union-bd}
    \sum_{\frac{m\e}{\q}<k\leq \frac{m\e \log m}{\q}} \P(N_k \geq k) \leq m \exp\left( -\log (2m/\d) \right) = \d/2.
\end{align}
When $k > \frac{m\e \log m}{\q}$, the upper bound (\ref{Nk-upper-bound}) gives 
\begin{align*}
    \E \widetilde{N}_k \leq k + \frac{\L m\e^2}{2}-\L\q k\e = k\left(1+\frac{\L m\e^2}{2k} - \L\q \e\right) \leq k\left(1 - \frac{\L\q \e}{2} \right),
\end{align*}
for $m$ large enough, since $\frac{\L m\e^2}{2k} \leq O\left(\frac{\e}{\log m}\right)$. Again using $\frac{1}{1-x} \geq 1+x$, this upper bound on $\E \widetilde{N}_k$ implies
\begin{align*}
    \P\left(\widetilde{N}_k \geq \E \widetilde{N}_k \cdot \frac{k}{\E \widetilde{N}_k} \right) &\leq \P\left(\widetilde{N}_k \geq \E \widetilde{N}_k \cdot \frac{k}{k\left(1-\frac{\L\q \e}{2} \right)} \right) \\
    &\leq \P\left(\widetilde{N}_k \geq \E \widetilde{N}_k \cdot \left( 1+\frac{\L\q \e}{2}\right) \right) \\
    &\leq \exp\left( -\frac{1}{3}\cdot \E \widetilde{N}_k \cdot \left(\frac{\L\q \e}{2} \right)^2 \right) \tag{by Lemma \ref{lemma-multiplicative-chernoff}}\\
    &\leq \exp\left( -\frac{1}{3}\cdot \frac{m\e}{2\q} \cdot \left(\frac{\L\q \e}{2} \right)^2 \right) \tag{by (\ref{Nk-lower-bound})} \\
    &= \exp\left( -\frac{\q \L^2 m\e^3}{24} \right).
\end{align*}
Since $\d \leq 1$, the above implies that the second piece of (\ref{switch-union-bound}) is bounded,
\begin{align*}
    \sum_{k > \frac{m\e \log m}{\q}} \P(N_k \geq k) \leq m \exp\left( - \log^3 (2m/\d) \right) \leq \d/2.
\end{align*}
Together with (\ref{small-k-union-bd}), we have shown
\begin{align*}
    \P(\hat{\tau}_\q > \tau^*_\q+\e) \leq \sum_{k > \frac{m\e}{\q}} \P(N_k \geq k) \leq \d.
\end{align*}
It remains to verify (\ref{Nk-upper-bound}) and (\ref{Nk-lower-bound}). To show (\ref{Nk-upper-bound}), note that for any $t \in [\tau^*_\q,\tau^*_\q+\e]$, the mean value theorem gives
\begin{align*}
\bar{f}(t) - \bar{f}(\tau^*_\q) \leq  - \L (t-\tau^*_\q)
\end{align*}
since $\bar{f}' \leq -\L$ on $[\tau^*_\q,\tau^*_\q+\e]$. Since $\bar{f}$ is decreasing on $[\tau^*_\q,\tau^*_\q+\q]$, this implies 
\begin{align*}
\bar{f}(t) \leq \begin{cases}
\bar{f}(\tau^*_\q) - \L (t-\tau^*_\q) \hspace{1em} &\tau^*_\q \leq t \leq \tau^*_\q+\e \\
\bar{f}(\tau^*_\q)-\L \e & \tau^*_\q+\e < t \leq \tau^*_\q+\q.
\end{cases}
\end{align*}
Thus the expectation can be bounded,
\begin{align*}
\E \widetilde{N}_k &= m \int_{\tau^*_\q}^{\tau^*_\q+\frac{\q k}{m}} \bar{f}(t) \de t \\
&= m\int_{\tau^*_\q}^{\tau^*_\q+\e} \bar{f}(t) \de t + m\int_{\tau^*_\q+\e}^{\tau^*_\q+\frac{\q k}{m}} \bar{f}(t) \de t \\
&\leq m\int_{\tau^*_\q}^{\tau^*_\q+\e} \bar{f}(\tau^*_\q) -\L(t-\tau^*_\q) \de t + m \int_{\tau^*_\q+\e}^{\tau^*_\q+\frac{\q k}{m}} \bar{f}(\tau^*_\q)-\L\e \de t \\
&= m\left[ \bar{f}(\tau^*_\q) \cdot \frac{\q k}{m} - \frac{\L(t-\tau^*_\q)^2}{2}  \bigg\vert^{\tau^*_\q+\e}_{\tau^*_\q} -  \L \e\left(\frac{\q k}{m}-\e \right) \right]\\
&= k - \frac{\L m\e^2}{2} - \L\q k \e + \L m \e^2 = k + \frac{\L m\e^2}{2} - \L\q k \e,
\end{align*}
which shows (\ref{Nk-upper-bound}). For (\ref{Nk-lower-bound}), note that the mean value theorem and the condition $\bar{f}' \geq -\L^{-1}$ on $[\tau^*_\q,\tau^*_\q+\e]$ imply that $\bar{f}(t) \geq \bar{f}(\tau^*_\q) - \L^{-1}(t-\tau^*_\q)$ for any $t \in [\tau^*_\q,\tau^*_\q+\e]$. Thus we have
\begin{align*}
\E \widetilde{N}_k &= m \int_{\tau^*_\q}^{\tau^*_\q+\frac{\q k}{m}} \bar{f}(t) \de t \\
&\geq m \int_{\tau^*_\q}^{\tau^*_\q+\e} \left( \bar{f}(\tau^*_\q)- \L^{-1}(t-\tau^*_\q)\right)\de t \\
&= m \e \bar{f}(\tau^*_\q) - \frac{m\e^2}{2\L} \\
&= \frac{m\e}{\q} - \frac{m\e^2}{2 \L} \geq \frac{m\e}{2\q},
\end{align*}
since for $m$ larger than some constant $C(\q,\L,\d)>0$, we have $\frac{m}{\log^3(2m/\d)} \geq 24 \q^2/\L^5$, which is equivalent to the last inequality above.
\end{proof}

A high probability lower bound can be shown under an extended monotonicity constraint of $f$ over the interval $(0,\tau^*_\q)$, as described in the next lemma.

\begin{lemma}
\label{lemma-threshold-concentration-lower}
Let $\d>0$. Suppose $f$ is decreasing on the interval $(0,\tau^*_\q)$ and that there exists some $\L>0$ for which $|f'(t)| \geq \L$ for all $t$ with $|t-\tau^*_\q| \leq \e$, where $\e \coloneqq \left(\frac{48}{\q \L^2}\right)^{1/3} m^{-1/3}\log(2m/\d)$. Then
\begin{align*}
\P(\hat{\tau}_\q < \tau^*_\q-\e) \leq \d, 
\end{align*}
for any $m \geq C(\q,\L,\d)$, a constant depending only on $\q,\L$ and $\d$.
\end{lemma}

\begin{proof}

Define $i^*$ as in Lemma \ref{lemma-switch}. If $\hat{\tau}_\q < \tau^*_\q -\e$, then there exists some $0\leq k \leq i^*$ for which $\hat{\tau}_\q = p_{(i^*-k)}$ and thus
\begin{align*}
    p_{(i^*-k)}-\frac{\q (i^*-k)}{m} \leq p_{(i^*)}-\frac{\q i^*}{m} \hspace{1em} \text{and} \hspace{1em} p_{(i^*-k)} < \tau^*_\q -\e .
\end{align*}
Since $p_{(i^*)}\leq \tau^*_\q$, it follows that the probability can be bounded,
\begin{align}
\nonumber
    \P(\hat{\tau}_\q < \tau^*_\q - \e) &\leq \P \left( \bigcup_{k=0}^m \left\{ p_{(i^*-k)} \leq \left(\tau^*_\q -\frac{\q k}{m}\right) \wedge (\tau^*_\q-\e) \right\} \cap\{i^* \geq k\} \right) \\
    \label{union-bd-k-small}
    &\leq \P \left(\bigcup_{0\leq k \leq \frac{m\e}{\q}}\left\{p_{(i^*-k)} \leq \tau^*_\q - \e \right\}\cap\{i^*\geq k\} \right) \\
    \label{union-bd-k-large}
    &+ \P \left(\bigcup_{k > \frac{m\e}{\q}}\left\{p_{(i^*-k)} \leq \tau^*_\q - \frac{\q k}{m}\right\}\cap \{i^*\geq k\} \right).
\end{align}
For (\ref{union-bd-k-small}), note that 
\begin{align*}
    p_{(i^*-k)} \leq \tau^*_\q - \e \Rightarrow N_\e\coloneqq \sum_{j=1}^m 1_{\{p_j \in [\tau^*_\q-\e,\tau^*_\q]\}} \leq k,
\end{align*}
since if at least $i^*-k$ of the $p$-values fall below $\tau^*_\q-\e$, and exactly $i^*$ of the $p$-values are below $\tau^*_\q$, then at most $k$ of the $p$-values fall in the interval $[\tau^*_\q-\e,\tau^*_\q]$. Since the $p$-values are independent, we again have $N_\e \sim \text{Generalized-Binomial}$ with sample size $m$ and average success probability $\bar{F}(\tau^*_\q)-\bar{F}(\tau^*_\q-\e)$. By the mean value theorem, for some $\xi \in [\tau^*_\q-\e,\tau^*_\q]$, we have
\begin{align*}
    \E N_\e = m(\bar{F}(\tau^*_\q)-\bar{F}(\tau^*_\q-\e)) = m\bar{f}(\xi) \e \geq m \bar{f}(\tau^*_\q)\e \geq k,
\end{align*}
since $\bar{f}$ is decreasing on $(0,\tau^*_\q)$, $\bar{f}(\tau^*_\q) = \q^{-1}$, and $k \leq \frac{m\e}{\q}$. It follows from Theorem 5 in \cite{hoeffding1956distribution} that
\begin{align}
\label{hoeffding-lower-tail}
    \P(p_{(i^*-k)}\leq \tau^*_\q -\e , i^*\geq k ) \leq \P(N_\e \leq k) \leq \P(\widetilde{N}_\e \leq k ),
\end{align}
where $\widetilde{N}_\e \sim \text{Binomial}(m,\bar{F}(\tau^*_\q)-\bar{F}(\tau^*_\q-\e))$. Further note that for any $t \in [\tau^*_\q-\e,\tau^*_\q]$, the mean value theorem and the condition $\bar{f}' \leq -\L$ on $[\tau^*_\q-\e,\tau^*_\q]$ imply
\begin{align*}
    \bar{f}(\tau^*_\q)-\bar{f}(t) = \bar{f}'(\xi) (\tau^*_\q-t) \leq -\L(\tau^*_\q-t),
\end{align*}
which further implies the following lower bound on the mean,
\begin{align}
\nonumber
    \E \widetilde{N}_\e &= m \int_{\tau^*_\q-\e}^{\tau^*_\q} \bar{f}(t) \de t \\
    \nonumber
    &\geq m \int_{\tau^*_\q-\e}^{\tau^*_\q} \bar{f}(\tau^*_\q) + \L(\tau^*_\q-t) \de t \\
    \label{Ne-lower}
    &= m \bar{f}(\tau^*_\q) \e - \frac{m\L}{2}(\tau^*_\q-t)^2\bigg\vert^{\tau^*_\q}_{\tau^*_\q-\e} = \frac{m\e}{\q} + \frac{m\L\e^2}{2}.
\end{align}
It follows that (\ref{hoeffding-lower-tail}) is bounded,
\begin{align*}
    \P(\widetilde{N}_\e \leq k) &= \P\left(\widetilde{N}_\e \leq \E\widetilde{N}_\e \cdot \frac{k}{\E\widetilde{N}_\e}\right) \\
    &\leq \P\left(\widetilde{N}_\e \leq \E\widetilde{N}_\e \cdot \frac{k}{\frac{m\e}{\q}\left(1+\frac{\L\q \e}{2}\right)}\right) \\ 
    &\leq \P\left(\widetilde{N}_\e \leq \E\widetilde{N}_\e \cdot \frac{1}{1+\frac{\L\q \e}{2}}\right) \tag{$k\leq \frac{m\e}{\q}$}.
\end{align*}
Now since $\frac{1}{1+x} \leq 1-x/2$ for $x \in [0,1]$, and since $\frac{\L\q \e}{2} \leq 1$ for $m$ larger than a constant, the above is bounded
\begin{align*}
    &\leq \P\left(\widetilde{N}_\e \leq \E \widetilde{N}_\e \left( 1-\frac{\L\q \e}{4}\right) \right) \\
    &\leq \exp\left( -\frac{1}{3}\cdot \E \widetilde{N}_\e \cdot \left(\frac{\L\q\e}{4} \right)^2\right) \tag{Lemma \ref{lemma-multiplicative-chernoff}} \\
    &\leq \exp\left(-\frac{1}{3} \cdot \frac{m\e}{\q} \cdot \left(\frac{\L\q\e}{4} \right)^2\right),
\end{align*}
since (\ref{Ne-lower}) implies $\E \widetilde{N}_\e \geq \frac{m\e}{\q}$. Plugging the definition of $\e$, we have shown
\begin{align*}
    \P(\widetilde{N}_\e \leq k) \leq \exp\left(-\frac{\q \L^2 m\e^3}{48} \right) = \exp\left(-\log^3(2m/\d) \right)\leq \frac{\d}{2m},
\end{align*}
so by the union bound, (\ref{union-bd-k-small}) is no larger than $\d/2$. 

For (\ref{union-bd-k-large}), similar to the first step in the analysis of (\ref{union-bd-k-small}), we have the implication
\begin{align*}
    p_{(i^*-k)}\leq \tau^*_\q - \frac{\q k}{m} \Rightarrow N_k \coloneqq \sum_{j=1}^m 1_{\{p_j \in [\tau^*_\q-\frac{\q k}{m},\tau^*_\q]\}} \leq k.
\end{align*}
We have $N_k \sim \text{Generalized-Binomial}$ with sample size $m$ and average success probability $\bar{F}\left(\tau^*_\q)-\bar{F}(\tau^*_\q - \frac{\q k}{m}\right)$ because the $p$-values are independent. By the mean value theorem, for some $\xi \in [\tau^*_\q - \frac{\q k}{m}, \tau^*_\q]$, we have
\begin{align*}
    \E N_k = m \bar{f}(\xi) \cdot \frac{\q k}{m} \geq k,
\end{align*}
since $\bar{f}$ is decreasing on $(0,\tau^*_\q)$ and $\bar{f}(\tau^*_\q)=\q^{-1}$. It thus follows from Theorem 5 in \cite{hoeffding1956distribution} that
\begin{align*}
    \P\left(p_{(i^*-k)}\leq \tau^*_\q - \frac{\q k}{m}, i^* \geq k\right) \leq \P(N_k \leq k) \leq \P(\widetilde{N}_k \leq k), 
\end{align*}
where $\widetilde{N}_k \sim \text{Binomial}\left(m,\bar{F}(\tau^*_\q)-\bar{F}\left(\tau^*_\q - \frac{\q k}{m} \right) \right)$.
For any $t \in [\tau^*_\q-\e,\tau^*_\q]$, the mean value theorem gives
\begin{align*}
\bar{f}(\tau^*_\q)-\bar{f}(t) \leq  - \L (\tau^*_\q-t)
\end{align*}
since $\bar{f}' \leq -\L$ on $[\tau^*_\q-\e,\tau^*_\q]$. Since $\bar{f}$ is decreasing on $(0,\tau^*_\q)$, this implies 
\begin{align*}
\bar{f}(t) \geq \begin{cases}
\bar{f}(\tau^*_\q) + \L (\tau^*_\q-t) \hspace{1em} &\tau^*_\q-\e\leq t \leq \tau^*_\q \\
\bar{f}(\tau^*_\q) + \L \e &t < \tau^*_\q-\e.
\end{cases}
\end{align*}
Thus $\E \widetilde{N}_k$ is bounded below,
\begin{align*}
    \E \widetilde{N}_k &= m \int_{\tau^*_\q - \frac{\q k}{m}}^{\tau^*_\q} \bar{f}(t) \de t \\
    &= m \int_{\tau^*_\q - \frac{\q k}{m}}^{\tau^*_\q-\e} \bar{f}(t) \de t + m \int_{\tau^*_\q -\e}^{\tau^*_\q} \bar{f}(t) \de t \tag{$k>\frac{m\e}{\q}$}\\ 
    &\geq m \int_{\tau^*_\q - \frac{\q k}{m}}^{\tau^*_\q-\e} (\bar{f}(\tau^*_\q)+\L \e) \de t + m \int_{\tau^*_\q -\e}^{\tau^*_\q} f(\tau^*_\q)+\L(\tau^*_\q-t) \de t \\
    &= m \bar{f}(\tau^*_\q) \cdot \frac{\q k}{m} + \L m\e \left(\frac{\q k}{m}-\e \right) - \frac{m\L}{2}(\tau^*_\q-t)^2\bigg\vert^{\tau^*_\q}_{\tau^*_\q -\e} \\
    &= k + \L\q k \e - m\L \e^2 + \frac{m\L \e^2}{2}.
\end{align*}
Simplifying, we have shown
\begin{align}
    \nonumber
    \E \widetilde{N}_k &\geq k + \L\q k \e - \frac{m\L \e^2}{2} \\
    \nonumber
    &> k + \L\q k \e - \frac{\L \q k\e}{2} \tag{$m\e < \q k$} \\
    \label{Nk-lower}
    &= k \left(1+\frac{\L \q \e}{2} \right).
\end{align}
Now since $\frac{1}{1+x}\leq 1-x/2$ for $x \in [0,1]$, and since $\frac{\L\q\e}{2} \leq 1$ for $m$ larger than a constant, we have
\begin{align*}
    \P(\widetilde{N}_k \leq k) &= \P\left(\widetilde{N}_k \leq \E \widetilde{N}_k \cdot \frac{k}{\E \widetilde{N}_k} \right) \\
    &\leq \P\left(\widetilde{N}_k \leq \E \widetilde{N}_k \cdot \left(1 - \frac{\L\q \e}{4} \right) \right) \\
    &\leq \exp\left( -\frac{1}{3}\cdot \E \widetilde{N}_k \cdot \left(\frac{\L\q\e}{4} \right)^2\right) \\
    &\leq \exp\left(-\frac{1}{3} \cdot \frac{m\e}{\q} \cdot \frac{\L^2\q^2 \e^2}{16} \right),
\end{align*}
since (\ref{Nk-lower}) together with $k>\frac{m\e}{\q}$ imply $\E \widetilde{N}_k \geq \frac{m\e}{\q}$. Plugging in the definition of $\e$, we have shown
\begin{align*}
    \P(\widetilde{N}_\e \leq k) \leq \exp\left(-\frac{\q \L^2 m\e^3}{48} \right) = \exp\left(-\log^3(2m/\d) \right)\leq \frac{\d}{2m},
\end{align*}
so by the union bound, (\ref{union-bd-k-large}) is no larger than $\d/2$. Since we've now shown that both terms (\ref{union-bd-k-small}) and (\ref{union-bd-k-large}) are below $\d/2$, the proof is complete.
\end{proof}

\begin{lemma}
\label{lemma-multiplicative-chernoff}
Let $X \sim \text{Binomial}(n,p)$. Then for any $0<\d<1/2$, we have
\begin{align*}
    \P(X \geq np(1+\d)) &\leq \exp\left(-\frac{1}{3} np\d^2 \right).
\end{align*}
\end{lemma}

\begin{proof}
By Markov's inequality, for any $t\geq 0$ we have
\begin{align*}
    \P(X \geq np(1+\d)) \leq \frac{\E e^{tX}}{e^{tnp(1+\d)}} = \frac{(1-p+pe^t)^n}{e^{tnp(1+\d)}} \leq \exp\left(np(e^t-1)-tnp(1+\d) \right).
\end{align*}
Letting $t=\log(1+\d)$, we have
\begin{align*}
    \P(X \geq np(1+\d)) \leq e^{np(\d-(1+\d)\log(1+\d))}.
\end{align*}
Now since $(1+\d)\log(1+\d)\geq \d+\frac{1}{3}\d^2$ for any $\d \in (0,1/2)$, we obtain the result. 
\end{proof}

\bibliographystyle{dcu}
\bibliography{reference.bib} 

\end{document}